\documentclass[12pt, reqno, twoside]{amsart}
\usepackage{amsmath, amsthm, amscd, amsfonts, amssymb, graphicx, color}
\usepackage[bookmarksnumbered, colorlinks, plainpages]{hyperref}
\hypersetup{colorlinks=true,linkcolor=red, anchorcolor=green, citecolor=cyan, urlcolor=red, filecolor=magenta, pdftoolbar=true}
\usepackage{enumitem}
\newtheorem{theorem}{Theorem}[section]
\newtheorem{lemma}[theorem]{Lemma}
\newtheorem{corollary}[theorem]{Corollary}
\theoremstyle{definition}

\newtheorem{proposition}[theorem]{Proposition}

\theoremstyle{remark}
\newtheorem{remark}[theorem]{Remark}
\numberwithin{equation}{section}
 \usepackage{mathrsfs}
\begin{document}

\setcounter{page}{1}

\title[Gradient estimates for nonlinear heat equations with compact boundary]{Some gradient estimates for nonlinear heat-type equations  on  smooth metric measure spaces with compact boundary}

\author[]{Abimbola Abolarinwa}

\address{Mathematical Analysis and Applications Research Group (MAARG), Department of Mathematics,
University of Lagos,
Nigeria.}
\email{\textcolor[rgb]{0.00,0.00,0.84}{a.abolarinwa1@gmail.com, aabolarinwa@unilag.edu.ng}}

 

 \date{July 25,  2023}

\begin{abstract}
In this paper we prove some Hamilton type and Li-Yau type gradient estimates on positive solutions to generalized nonlinear parabolic equations on smooth metric measure space with compact boundary.
The geometry of the space in terms of lower bounds on the weighted Bakry-\'Emery Ricci curvature tensor and weighted mean curvature of the boundary are key in proving generalized local and global gradient estimates.  Various applications of these gradient estimates in terms of parabolic Harnack inequalities and Liouville type results are discussed.  Further consequences to some specific models informed by the nature of the nonlinearities are highlighted.
\end{abstract}

 \maketitle
 
\tableofcontents

\section{Introduction}

Let the triple $(M,g,e^{-\phi}dv)$ be a smooth metric measure space whereby the pair $(M,g)$ is a complete Riemannian manifold of dimension greater than or equal to $2$ equipped with volume measure $dv$ and $\phi \in \mathscr{C}^\infty(M)$ is a smooth potential function. The smooth metric measure space $M:= (M,g,e^{-\phi}dv)$ is linked with a second order self-adjoint diffusion operator called $\phi$-Laplacian, $\mathscr{L}_\phi$,  which is usually associated with weighted Ricci tensor, $\mathscr{R}ic_\phi$.  In this paper we are interested in gradient estimates of parabolic and elliptic types and their consequences for positive solutions to a class of nonlinear heat-type equations on $M$ with compact boundary $\partial M$.  Specifically,  our attention is on the positive smooth solutions $u=u(x,t)$ with $x\in M$ and $t\in (-\infty, +\infty)$ satisfying the nonlinear parabolic partial differential equation
\begin{align}\label{11}
\begin{cases}
\displaystyle \frac{\partial u}{\partial t}&= \mathscr{L}_\phi u +H(u)\\
\displaystyle u |_{\partial M} & = \text{const}.
\end{cases}
\end{align}
Here $H=H(s)$ is a sufficiently smooth nonlinear function for $s>0$.
Our interest is motivated by the  classical gradient estimates pioneered by \cite{[Ha93],[LY86]} which have been of great importance in the  analysis of partial differential equations and geometric analysis.   The antecedence of gradient estimates and their accompanied consequences in the field of geometric analysis and partial differential equations are well studied, and they can not be underestimated.   Gradient estimates have  served as a basic tool in the study of linear and nonlinear PDEs  yielding several applications ranging from Harnack estimates to Liouville theorems, Heat kernel estimates to eigenvalue bounds,  existence, compactness and regularity of solutions,  and so on.  See  \cite{[Ab20],[Ab19],[Ab1],[Ab2],[Ab0],[Ba1],[CZh18],[Li],[MZ],[SZ06],[Wu15],[Ya]}  and the references cited therein for some related theories and applications.  

Specifically, Li and Yau \cite{[LY86]} established local and global version of gradient type Harnack inequality for positive solutions to the heat equation
\begin{align}\label{es1}
\frac{\partial u}{\partial t} =  \Delta u
\end{align}
on complete Riemannian manifolds with or without boundary.  Their inequality was used to derive bounds for the fundamental solution  of \eqref{es1} under the nonnegativity assumption on the Ricci curvature tensor.   Hamilton in \cite{[Ha93]} proved another version of Li-Yau estimate which is time independent for the heat equation \eqref{es1} on compact Riemannian manifold.  Hamilton's result allows to compare temperatures of two different points at the same time.  Inspired by Hamilton's result, Souplet and Zhang \cite{[SZ06]} got different version of gradient estimate which is also time independent yielding a new Liouville-type theorems of Cheng-Yau type \cite{[CY75]} for the heat equation on non-compact Riemannian manifolds. Souplet-Zhang's result allows the comparison of temperature distribution simultaneously, without any lag in time, even for noncompact manifolds.
Since the appearance of these works by Li and Yau \cite{[LY86]}, Hamilton \cite{[Ha93]} and Souplet and Zhang \cite{[SZ06]} several results on generalizing, extending and improving  classical gradient estimates with several applications have appeared in literature.,  including but not limited to those ones cited above.  See also \cite{[Ab22],[AE],[CCK],[CLPW],[CaM1],[CaM2],[CY75],[Du23],[Du],[DKA],[GH],[Hou],[Li05],[Ma],[Po1],[Po2],[SZ06],[Wu],[Ya1],[Ya2]}.  The gradient estimates discussed in this paper are broadly categorized into two, namely, elliptic type and parabolic type. The elliptic type is also known as Hamilton-Souplet-Zhang type, they only involve the solution and its gradient. The parabolic type which involve the time derivative of the solution in addition is known as Li-Yau type.

The study of evolution equations and their steady states involving the weighted Laplacian on manifolds on its own has been of tremendous usefulness in recent research development due to their profound intimacy and connection with differential geometry, probability theory, quantum theory and statistical mechanics to list a few.  An important aspect of this research line is the linkage between the diffusion operator $\mathscr{L}_\phi$ and the geometry of the underlying manifold in the form of generalized weighted Ricci curvature tensor, popularly known as Bakry-\'Emery-Ricci tensor $\mathscr{R}ic_\phi^m$ in the literature \cite{[BGL],[Lot],[Vil]}. This enables convenient and nice generalizations of manifold results such as Li-Yau Harnack inequalities,  Laplacian comparison,  Myers' compactness,  Bishop-Gromov’s volume comparison and Cheeger-Gromoll's splitting theorems, e.t.c.  to smooth metric measure spaces \cite{[BE],[BGL],[Bri],[Li05],[LLi1],[LLi2],[LLi3],[Lot],[Sak2],[Vil],[WZZ16],[WeW09]}. The generalizations are not usually trivial, as  they  may involve a careful consideration of properties of the potential function and the choice of input data such that special conditions would not be required on the potential function. 

The nonlinear evolution equation \eqref{11} has attracted attention of analysts and applied mathematicians owing to its numerous applications in geometry, mathematical physics and several other areas of science. The simplest case being when $H(u)\equiv 0$ in which case \eqref{11} is just the weighted heat equation, $(\partial_t - \mathscr{L}_\phi)u=0$,  which is archetypical diffusion processes.  For instance it is linked with Markov processes, and the fundamental solution of the weighted heat equation is actually the transition probability density function of the diffusion processes associated with $\mathscr{L}_\phi$ which can be constructed using It\^o calculus via certain stochastic differential equation \cite{[BGL],[Li05],[LLi3]}.
When we have $H(u)=au(1-u)$ or $H(u)=au(1-u^2)$ for $a>0$ and $0<u\le 1$,  \eqref{11} is known in literature as Fisher-KPP or Allen-Cahn equation, both of which arise in biological and chemical dynamics \cite{[AM20],[AC],[CC],[XChen],[Fi],[KPP]}. When $H(u)=au\log u$ \eqref{11} is closely related to the gradient Ricci solitons and weighted Log-Sobolev inequalities on smooth metric measure spaces.  The stationary form  of \eqref{11} is related to weighted Yamabe-type equation for $H(u)=\lambda u+\mu u^p$ which is applicable in conformal geometry,  and also Einstein Scalar field equation for $H(u)= au+bu^\alpha + cu^\beta$ which is widely applicable in Einstein field theory of general relativity \cite{[Cas1],[Cas2],[DKA]}.  In general, the nonlinear parabolic equation \eqref{11} is the so-called weighted reaction-diffusion equation which can be found in many mathematical models in physics, biology, chemistry, engineering in which case the term $\mathscr{L}_\phi$ accounts for diffusion while the term $H(u)$ is the reaction term.

Motivated by the above works, most especially, Li-Yau \cite{[LY86]}, Hamilton \cite{[Ha93]}, Souplet-Zhang \cite{[SZ06]},  Li \cite{[Li]}, Fu and Wu \cite{[FW21]} and Abolarinwa, Ali and Mofarreh \cite{[Ab22]},  we prove local  and global elliptic-type and parabolic-type gradient estimates for positive solutions to generalized nonlinear parabolic and elliptic partial differential equations on smooth metric measure space with boundary satisfying some lower bound assumptions on $m$-Bakry-\'Emery Ricci tensor and weighted mean curvature of the boundary. Applications of these estimates are made towards deriving the Harnack inequalities and Liouville-type results.

This introductory section is set to provide an overview of the theme of and motivation for this study. Previous studies and background information on the topic are highlighted. Meanwhile, the rest of this paper is structured as follows:

\begin{itemize}
\item Section \ref{sec2}: {\bf Preliminaries and statement of main results}: The purpose of this section is in two folds. First, the preliminaries subsection (Subsection \ref{sec21}) presents basic concepts and fundamental results of smooth metric measure spaces with boundary which are relevant to this study, as well as fixing the required notation.  More concretely, descriptions of weighted Laplacian, Bakry-\'Emery Ricci tensor, weighted Bochner formula, weighted Reilly-type formula, mean curvature of the boundary, Laplacian comparison result on smooth metric measure space with boundary are made. Secondly, the statement of main results subsection (Subsection \ref{sec22}) states the three main results on gradient estimates. The first two results (Theorem \ref{thm21a} and Theorem \ref{thm22}) are elliptic in nature and are called Hamilton-Souplet-Zhang  type gradient estimates, while the third result (Theorem \ref{thm23}) is on Li-Yau gradient estimates.

\item Section \ref{sec3}: {\bf Hamilton-Souplet-Zhang  type gradient estimates I}: Theorem \ref{thm21a} is proved in this section. Some fundamental lemmas needed in the proof are stated and their detailed proofs are given.

\item Section \ref{sec4}: {\bf Hamilton-Souplet-Zhang  type gradient estimates II}: Here in this section we prove Theorem \ref{thm22}.  We also state and prove some fundamental lemmas needed in the proof of Theorem \ref{thm22}. The approach adopted in this section is similar to that of Section \ref{sec3}, but with a different auxilliary function.

\item Section \ref{sec5}: {\bf Li-Yau type gradient estimates}: Theorem \ref{thm23} is proved in this section. The approach introduced by Li \cite{[Li]} is adopted by choosing a cut-off function different from the one applied in Sections \ref{sec3} and \ref{sec4}.

\item Section \ref{sec6}: {\bf Applications of Li-Yau  type gradient estimates}: This section presents some applications of Theorem \ref{thm23} to deriving global gradient estimates and Li-Yau (parabolic) Harnack inequalities for positive solutions of \eqref{11} as presented in Theorem \ref{thm62} and Theorem \ref{cor63}. Li-Yau type gradient estimates for weighted nonlinear elliptic equations on weighted manifold with boundary and Liouville-type theorems are derived as applications of the parabolic gradient etimstes. These results are consequently applied to special cases in the form of logarithmic and power-like nonlinearities (Subsection \ref{sec64}).

\item Section \ref{sec7}: {\bf Applications of elliptic gradient estimates}:  Similar to Section \ref{sec6}, this final section presents some applications of Hamilton-Souplet-Zhang gradient estimates. The applications yield global gradient estimates and Liouville-type results,  among others. These results are also applied to special cases in the form of logarithmic and power-like nonlinearities.
\end{itemize}

\section{Preliminaries and statement of main results}\label{sec2}

\subsection{Preliminaries}\label{sec21}
As stated earlier we shall consider  several gradient estimates for positive solutions of \eqref{11} on smooth metric measure spaces with compact boundary. Here in this section we want to give some preliminaries and notation as a way of formal introduction of a weighted manifold with compact boundary. 
Let $(M,g)$ be an $n$-dimensional  complete Riemannian manifold equipped with metric $g$ and its induced volume form $dv$. Let $\phi$ be a real valued $\mathscr{C}^\infty$-function on $M$ then the triple $(M,g,e^{-\phi}dv)$  endowed with weighted volume form $e^{-\phi}dv$ is called a smooth metric measure space (also called weighted manifold or manifold with density). The concept of smooth metric measure spaces was introduced by Bakry and \'Emery \cite{[BE]} in the mid 80's as a natural extension of smooth manifolds and has since been applied in several contexts such as  functional inequalities, optimal transport, Ricci solitons to mention but a few, see \cite{[BE],[BGL],[LV],[Pe02],[Vil]}.  

\subsubsection{Weighted Laplacian and weighted Ricci tensor}
The generalized Laplacian and generalized Ricci curvature tensor are defined on a  weighted manifold respectively by  $\phi$-Laplacian (also called Weighted Laplacian or Witten Laplacian) and Bakry-\'Emery Ricci curvature tensor $\mathscr{R}ic_\phi$.  For a $\mathscr{C}^2$ function $u$, the weighted Laplacian $\mathscr{L}_\phi$, which is a symmetric diffusion operator with respect to the invariant volume form $e^{-\phi}dv$, is given by
$$\mathscr{L}_\phi u := e^{\phi}\text{div}(e^{-\phi} \nabla u) = \Delta u - \langle \nabla \phi, \nabla u \rangle,$$
while the weighted $m$-Bakry-\'Emery Ricci curvature tensor $\mathscr{R}ic_\phi^m$ is given by 
\begin{align}\label{e0}
\mathscr{R}ic_\phi^m & =  \mathscr{R}ic + \text{Hess} \phi  - \frac{\nabla \phi \otimes \nabla \phi}{m-n}, \ \ \ m \in [n, \infty)
\end{align}  
Here, 
$$\text{div} (X) = \frac{1}{\sqrt{|g|}}\sum_{i=1}^n \frac{\partial}{\partial x^i}(\sqrt{|g|}X^i),  \ \ \ \ (\nabla u)^i = \sum_{i=1}^n g^{ij} \frac{\partial u}{\partial x^j} $$
and
$$\Delta u = \text{div}(\nabla u) =  \frac{1}{\sqrt{|g|}}\sum_{i,j=1}^n \frac{\partial}{\partial x^i} \left(\sqrt{|g|} g^{ij} \frac{\partial u}{\partial x^j} \right),$$
are respectively the usual divergence of a smooth vector field $X$,  gradient and Laplace-Beltrami operators on $(M,g)$ (in local coordinates),  where $|g|=\det(g_{ij})$ and $[g^{ij}] = [g_{ij}]^{-1}$ denote the determinant and the inverse of the component matrix of the metric tensor $g$. While $\langle \cdot, \cdot \rangle$, $\text{Hess} \phi=\nabla^2\phi$ and $\mathscr{R}ic$ are respectively the inner product with respect to  $g$, Hessian of function $\phi$ and the Ricci curvature tensor of $M$. The case $m=n$ is admissible in formula \eqref{e0} only when $\phi$ is a constant function.  The Bakry-\'Emery Ricci  tensor is interpreted as the limit of $ \mathscr{R}ic_\phi^m$ as $m\to \infty$  and one simply writes 
$$\lim_{m\to \infty} \mathscr{R}ic_\phi^m = \mathscr{R}ic + \text{Hess} \phi=:\mathscr{R}ic_\phi.$$

\subsubsection{Weighted Bochner formula}
It is important to note that $\phi$-Laplacian and the Bakry-\'Emery Ricci  tensor are connected by the weighted Bochner formula
\begin{align}\label{e1}
\mathscr{L}_\phi(|\nabla v|^2) = 2\|\text{Hess}\ v\|_{HS}^2 + 2\langle \nabla v, \nabla \mathscr{L}_\phi v \rangle + 2 \mathscr{R}ic_\phi(\nabla v, \nabla v),
\end{align}
where $v$ is  $\mathscr{C}^\infty$ and $\|\cdot \|_{HS}$ is Hilbert Schmidt norm with respect to   metric $g$.
Applying the elementary inequality of the form 
$(a+b)^2\ge \frac{a^2}{1+s} - \frac{b^2}{s}, \ s>0$ 
and a Newton inequality in the form $\|\text{Hess}\ v\|^2 \ge \frac{1}{n}(\Delta v)^2$, we clearly have
\begin{align}\label{ee1}
\|\text{Hess}\ v\|_{HS}^2 +  \mathscr{R}ic_\phi(\nabla v, \nabla v) &\ge \frac{1}{n}(\mathscr{L}_\phi v +\langle\nabla \phi, \nabla v\rangle)^2 + \mathscr{R}ic_\phi(\nabla v, \nabla v)\nonumber \\
&\ge \frac{1}{n}\left(\frac{n(\mathscr{L}_\phi v)^2}{m} -\frac{n \langle\nabla \phi, \nabla v\rangle^2}{m-n}   \right) + \mathscr{R}ic_\phi(\nabla v, \nabla v)\nonumber \\
&= \frac{1}{m}(\mathscr{L}_\phi v)^2 + \mathscr{R}ic^m_\phi(\nabla v, \nabla v).
\end{align}
Then \eqref{e1} yields a very important inequality that would be applied under the name weighted Bochner inequality
\begin{align}\label{e1a}
\mathscr{L}_\phi(|\nabla v|^2) \ge \frac{2}{m}(\mathscr{L}_\phi v)^2  + 2\langle \nabla v, \nabla \mathscr{L}_\phi v \rangle + 2 \mathscr{R}ic^m_\phi(\nabla v, \nabla v).
\end{align}
The identity \eqref{e1a} is equivalent to the curvature dimension condition $CD(\mathsf{k},m)$ with respect to the operator $\mathscr{L}_\phi$ under the lower bound assumption $\mathscr{R}ic^m_\phi \ge \mathsf{k}g$, $\mathsf{k}\in \mathbb{R}$, $m \in [n,\infty]$ embedded in the next given inequality
\begin{align}\label{e2a}
\frac{1}{2}\mathscr{L}_\phi(|\nabla v|^2)- \langle \nabla v, \nabla \mathscr{L}_\phi v \rangle \ge \frac{1}{m}(\mathscr{L}_\phi v)^2  +  \mathsf{k}|\nabla v|^2.
\end{align}
Obviously, the condition $\mathscr{R}ic_\phi\ge \mathsf{k}g$ implies $CD(\mathsf{k},\infty)$ and one then sees that \eqref{e2a} implies
\begin{align}\label{e3a}
\frac{1}{2}\mathscr{L}_\phi(|\nabla v|^2)- \langle \nabla v, \nabla \mathscr{L}_\phi v \rangle \ge  \mathsf{k}|\nabla v|^2 
\end{align}
but \eqref{e3a} does not imply \eqref{e2a}. Thus, the lower bound condition on $\mathscr{R}ic_\phi^m$ implies corresponding condition on $\mathscr{R}ic_\phi$ but not the other way round. For more details on Bochner formula and curvature dimension condition see \cite{[BE],[BGL],[Li05],[WeW09]}.  We note that the lower bound assumption is imposed on the generalized weighted Bakry-\'Emery Ricci tensor $\mathscr{R}ic_\phi^m$, $m\in [n,\infty)$ in all the reults of this paper.

The relation $\mathscr{R}ic_\phi=\kappa g$ for some constant $\kappa$ is just the gradient Ricci solitons equation which is fundamental in the singularity analysis of the Hamilton-Perelman Ricci flow \cite{[Ha],[Pe02]} (interested readers can find \cite{[Cao]} for detail exposition on Ricci solitons), whilst the relation $\mathscr{R}ic_\phi^m=\kappa g$ corresponds to the quasi-Einstein equation which consists gradient Ricci solitons for $m=\infty$, and Einstein metrics  when $\phi$ is constant. It also relates to warped product Einstein metrics for $m$ being a positive integer.  Quasi-Einstein metrics have been studied by several authors see \cite{[Cas1],[Cas2]} for examples.

\subsubsection{Mean curvature of the boundary}
Let  $\nu$  be the outward pointing unit vector normal  to $\partial M$. The weighted mean curvature (also called $\phi$-mean curvature) of the boundary of $M$  denoted by $\mathscr{H}_\phi$  is defined by 
$$\mathscr{H}_\phi:= \mathscr{H}-\langle\nabla \phi, \nu\rangle,$$
where $\mathscr{H}$ is the mean curvature of the boundary of $(M,g)$ with respect to $\nu$.  If $\nu$ is a unit normal pointing inward, the weighted mean curvature will be 
 $\mathscr{H}_\phi:= \mathscr{H}+\langle\nabla \phi, \nu\rangle$. 
Let $\vartheta_\nu$  be the normal derivative of a smooth function $\vartheta$ on $\partial M$ and $\rm{II}_\phi$  be the second fundamental form of $\partial M$ with respect to $\nu$. The following  identity which will be henceforth referred to as Reilly-type formula holds for all smooth function $\vartheta$ on ${M}$ (see \cite[Proposition 2.3]{[FW21]} or \cite{[DDW21],[DW21]})
\begin{align}\label{Rei}
(|\nabla \vartheta|^2)_\nu & = 2\vartheta_\nu\left[\mathscr{L}_\phi\vartheta - \mathscr{L}_{\partial M,f}(\vartheta|_{\partial M}) - \mathscr{H}_\phi\vartheta_\nu \right] + 2\langle\nabla_{\partial M}(\vartheta|_{\partial M}), \nabla_{\partial M}(\vartheta_\nu) \rangle|_{g_{\partial M}} \nonumber\\
& \hspace{.5cm} -2 \rm{II}_\phi \langle\nabla_{\partial M}(\vartheta|_{\partial M}), \nabla_{\partial M}(\vartheta|_{\partial M}) \rangle|_{g_{\partial M}},
\end{align}
where $\mathscr{L}_{\partial M,\phi}:=\Delta_{\partial M} - \langle \nabla_{\partial M} \phi \nabla_{\partial M}\cdot\rangle$ and $\nabla_{\partial M}$  are respectively the $\phi$-Laplacian and covariant derivative with respect to the induced metric $g_{\partial M}$ on  $\partial M$.
The formula \eqref{Rei} can be derived from the proof of the weighted Reilly formula in \cite{[MD10]}  (see also \cite[Appendix 4]{[CMJ]} while its classical version has been proved in \cite{[Rei]}.  One can also observe that \eqref{Rei}  does not involve Ricci curvature tensor which is due to the fact that weighted Ricci tensor bound concerns only the interior.  Note that only the mean curvature bound is imposed on the boundary.

\subsubsection{Laplacian comparison on weighted manifold with boundary}
Let 
$$ \varrho(x):= \varrho_{\partial M}(x)= d(x, \partial M)$$
be the distance function from the boundary $\partial M$ to any point $x$ in $M$, we denote by 
$$B_R(\partial M):=\{x:\   \varrho_{\partial M}(x) <R\}$$
for $R>0$ the $R$-neigbourhood of $\partial M$.  Based on the result in \cite{[Sak2]} one may suppose that $\varrho$ is smooth away from the cut locus for the boundary, $\text{Cut}(\partial M)$. By the results in \cite{[WZZ16]} and \cite{[Sak2]} (see also \cite[Theorem 2.1]{[DDW21]} and \cite[Theorem 2.1]{[FW21]}) we have  weighted Laplacian comparison theorem for the distance function $\varrho$ on $M$ under some assumptions. 
\begin{theorem}\label{thmW}(\cite{[Sak2]}
Let $M=(M,g, e^{-\phi}dv)$ be an $n$-dimensional complete smooth metric measure space with compact boundary $\partial M$.  Suppose the weighted Ricci tensor of $M$ and the weighted mean curvature of $\partial M$ satisfy 
$\mathscr{R}ic_\phi \ge -(n-1)\kappa$ and $\mathscr{H}_\phi \ge -\ell$ for some constants $\kappa\ge 0$ and $\ell\in \mathbb{R}$.
Then
\begin{align}\label{eq52b}
\mathscr{L}_\phi \varrho_{\partial M}(x) \le (n-1)\kappa R +\ell
\end{align}
for all $x \in B_R(\partial M)$ outside the $\rm{Cut}(\partial M)$.
\end{theorem}
The above comparison theorem  for the generalized Laplacian in smooth metric measure spaces with boundary  was proved in \cite{[Sak2]}. The authors in \cite{[WZZ16]} had earlier proved a similar result on weighted manifold with boundary under some assumption. This comparison theorem holds for all $R>0$ in contrast to Wei-Wylie comparison theorem \cite{[WeW09]} on weighted manifold with boundary.

\subsubsection{Cut-off functions}
Some space-time localization argument will be required so as to prove the elliptic type gradient estimates.  As a routine we do this through the construction of suitable cut-off function as in the next lemma  based on the idea introduced by Souplet and Zhang \cite{[SZ06]} for manifold without boundary. This is the technique used by Li and Yau \cite{[LY86]}, see also  \cite{[Ab19],[Ab1]}. 
In the setting of manifolds with boundary we refer to  \cite{[DDW21]},\cite{[FW21]} and \cite{[KS21]} for similar applications.

Define a set called parabolic (space-time) cylinder
$$\mathcal{Q}_{R,T}:= B_R(\partial M)\times [T_0-T,T_0] \subset \overline{M}\times(-\infty, \infty)$$
for $R\ge 2$ and $T>0$.  For a fixed point in $M$ at time $t$, we set 
$$\psi(x,t)=\bar{\psi}( \varrho_{\partial M}(x),t)$$
to be a cut-off function depending on space and time. The existence and other useful properties of $\bar{\psi}( \varrho,t)$ in this context are stated in the following lemma.
\begin{lemma}\label{lem23}
Let $M=(M,g,e^{-\phi}dv)$ be an $n$-dimensional complete smooth metric measure space with compact boundary $\partial M$.  There exists a smooth non-increasing cut-off-function $\psi=\psi(\rho_{\partial M},t)$ supported in the closure of $\mathcal{Q}_{R,T}(\partial M)$, $T, R>0$, such that for
\begin{enumerate}
\item  $\psi(x,t) = \bar{\psi}(\varrho_{\partial M}(x),t)$ with $0\leq \bar{\psi}(\varrho_{\partial M}(x),t) \leq 1$ in $[0,R]\times[T_0-T,T_0]$ and $\bar{\psi}$ is supported in a subset of $[0,R]\times[T_0-T,T_0]$.
\item $\bar{\psi}(\varrho_{\partial M}(x),t)=1$ and $\displaystyle \frac{\partial\bar{\psi}(\varrho, t)}{\partial \varrho}=0$ in $[0,R/2]\times[T_0-T,T_0]$.
\item   $\displaystyle \left|\frac{\partial\bar{\psi}}{\partial \varrho}\right|\leq \frac{C_a\bar{\psi}^a}{R}$ and $\displaystyle \left|\frac{\partial^2 \bar{\psi}}{\partial \varrho^2}\right| \leq \frac{C_a\bar{\psi}^a}{R^2}$, where $0<a<1$ and $C_a>0$ is a constant depending on $a$, hold on $[0,\infty) \times [T_0-T,T_0]$.
\item The estimate $\displaystyle \left|\frac{\partial\bar{\psi}}{\partial t} \right| \leq \frac{C\bar{\psi}^{1/2}}{T}$ holds on $[0,\infty) \times [T_0-T,T_0]$ for some constant $C>0$ independent of $R$ and $T$.
\end{enumerate}
\end{lemma}

\subsubsection{Some elementary inequalities}
Some important elementary inequalities will be applied in estimating certain bounds during the proof of our results, namely Cauchy-Schwarz, Young's and Newton inequalities.  For clarity we state them in a simple form, though they will be applied with little modification. 
Cauchy-Schwarz inequality says that the product of any two real numbers is not bigger than the product of their absolute value, that is,  if $a,b\in \mathbb{R}$, it holds that
$ab\le |a||b|$, where $|\cdot|$ is the absolute value.  Young inequality states that for any nonnegative real numbers $A,B\ge 0$ and real numbers $r,s>1$ such that $1/s+1/r=1$, it holds that 
$$AB \le (1/s)A^s +(1/r)B^r,$$
with equality if $A^s=B^r$. Combining Cauchy-Schwarz and Young's inequality therefore we have for product of any two real numbers $a$ and $b$
$$ab\le |a||b| \le (1/s)|a|^s +(1/r)|b|^r.$$

A Newton inequality (as applied in this paper) states that 
\begin{align}\label{NI}
\|\text{Hess}(f)\|^2 \ge (1/n)(\Delta f)^2.
\end{align}
This is derived as a consequence of the generalized AM-GM inequality $\|\mathbb{A}\|^2 \ge (1/n)(\text{tr} \mathbb{A})^2$, where $\mathbb{A}$ is an $n\times n$ real symmetric matrix. For instance, take the $n$-dimensional vector $\mathbf{A}=(A_{ii}), i=1,\cdots n,$ and apply Cauchy-Schwarz inequality we can derive
$$\left(\sum_{i=1}^n A_{ii}\right)^2 \le n \sum_{i=1}^n  (A_{ii})^2 \le n \sum_{i,j=1}^n (A_{ij})^2.$$

\subsection{Statement of main results}\label{sec22}
In the next we state the main results of this paper while their detailed proofs and applications are postponed till  later sections.  First consider the following notation that will appear in the results and their proofs. For time dependent smooth function $u=u(x,t)$ we denote time (partial) derivative by $u_t=\partial_tu=\frac{\partial u}{\partial t}$. For a real valued quantity $[h]$ we denote its positive and negative parts by $[h]_+:= \max\{0,h\}$ and $[h]_-:=\min\{0,h\}$, respectively. Ordinary derivative is denoted by $'$, that is $H'(u)=\frac{d}{du} H(u)$ and $H''(u)=\frac{d^2}{du^2}H(u)$. Furthermore, we denote $\widehat{H}(u)=H(u)/u$ and $H_u=H'(u)-H(u)/u$. The notation $\log$ represents natural logarithm. The letters $C,C_1, C_2, \cdots, $ represent universal positive constants which may depend only on $m$ and $n$, with varying value from line to line.

\subsubsection{Hamilton-Souplet-Zhang type gradient estimates I}

\begin{theorem}\label{thm21a}
Let $M=(M,g,e^{-\phi}dv)$ be an $n$-dimensional complete smooth metric measure space with compact boundary satisfying the bounds
\begin{align*}
\mathscr{R}ic_\phi^m \ge -(m-1)\kappa \ \ \text{and}\ \ \mathscr{H}_\phi \ge - \ell
\end{align*}
for some constants $\kappa, \ell \ge 0$. Let $u=u(x,t)$ be a smooth positive solution to \eqref{11} (i.e.   $(\partial_t -\mathscr{L}_\phi) u =H(u)$ in $\mathcal{Q}_{R,T}(\partial M$). Suppose further that $u$ satisfies the Dirichlet boundary condition (i.e., $u(\cdot,t)|_{\partial M}$ is constant for each time $t \in (-\infty,+\infty)$),
\begin{align*}
u_\nu\ge 0\ \ \text{and}\ \ u_t \le H(u)
\end{align*}
over $\partial M \times [T_0-T, T_0]$. Then there exists a constant $C>0$ depending only on $n,m$, such that for all $(x,t)\in \mathcal{Q}_{R/2,T/2}(\partial M)$
\begin{align}
\frac{|\nabla u|}{1+\sqrt{u}} & \le C\Bigg\{ \left(\frac{1}{R} + \frac{1}{\sqrt{T}}+\sqrt{L}\right)\Big(\sup_{\mathcal{Q}_{R,T}(\partial M) }\{u\}\Big) \nonumber\\
& \hspace{2cm} + \Big(\sup_{\mathcal{Q}_{R,T}(\partial M) }\{\sqrt{\mathbb{M}_+[u]}\}\Big) \Big(\sup_{\mathcal{Q}_{R,T}(\partial M) }\{\sqrt{u}\}\Big) \Bigg\} ,
\end{align}
where  $L= (\kappa^2+\ell^4)^{1/2}$,
$$\mathbb{M}_+[u] = \left[\frac{u[(1/2+\sqrt{u})\widehat{H}(u) - (1+\sqrt{u})H'(u)]}{(1+\sqrt{u})^2} \right]_+ \ \ \text{and}\ \ \widehat{H}(u):=\frac{H(u)}{u}.$$  
\end{theorem}

\subsubsection{Hamilton-Souplet-Zhang type gradient estimates II}
\begin{theorem}\label{thm22}
Let $M=(M,g,e^{-\phi}dv)$ be an $n$-dimensional complete smooth metric measure space with compact boundary satisfying the bounds
\begin{align*}
\mathscr{R}ic_\phi^m \ge -(m-1)\kappa \ \ \text{and}\ \ \mathscr{H}_\phi \ge - \ell
\end{align*}
for some constants $\kappa, \ell \ge 0$. Let $u=u(x,t)$ be a smooth positive solution to \eqref{11}  in $\mathcal{Q}_{R,T}(\partial M)$. Suppose further that $u$ satisfies the Dirichlet boundary condition (i.e., $u(\cdot,t)|_{\partial M}$ is constant for each time $t \in (-\infty,+\infty)$),
\begin{align*}
u_\nu\ge 0\ \ \text{and}\ \ u_t \le H(u)
\end{align*}
over $\partial M \times [T_0-T, T_0]$. Then there exists a constant $C>0$ depending only on $n,m$, such that for all $(x,t)\in \mathcal{Q}_{R/2,T/2}(\partial M)$
\begin{align}
\frac{|\nabla u|}{\sqrt{u}} & \le C\Bigg\{\frac{1}{R} + \frac{1}{\sqrt{T}}+\sqrt{L}+ \Gamma^H_+[u] \Bigg\}\Big(\sup_{\mathcal{Q}_{R,T}(\partial M) }\{\sqrt{u}\}\Big) ,
\end{align}
where  $L= (\kappa^2+\ell^4)^{1/2}$,
$$\Gamma^H_+[u] =  \sup_{\mathcal{Q}_{R,T}(\partial M) }\left\{\left[ 2H'(u)-\widehat{H} \right]_+ \right\}^{1/2} \ \ \text{and}\ \ \widehat{H}(u):=\frac{H(u)}{u}.$$  
\end{theorem}

\subsubsection{Li-Yau type gradient estimates}

\begin{theorem}\label{thm23}
Let $M=(M,g,e^{-\phi}dv)$ be an $n$-dimensional complete smooth metric measure space with compact boundary satisfying the bounds
\begin{align*}
\mathscr{R}ic_\phi^m \ge -(m-1)k \ \ \text{and}\ \ \mathscr{H}_\phi \ge - l
\end{align*}
for some constants $k,l \ge 0$.  Let $u=u(x,t)$ be a smooth positive solution to \eqref{11} in $\mathcal{Q}_{2R,T}(\partial M)$. Suppose further that $u$ satisfies the Dirichlet boundary condition (i.e., $u(\cdot,t)|_{\partial M}$ is constant for each time $t \in (-\infty,+\infty)$),
\begin{align*}
u_\nu\ge 0\ \ \text{and}\ \ u_t \le H(u)
\end{align*}
over $\partial M \times [T_0-T, T_0]$.   Then have for all $(x,t)\in \mathcal{Q}_{R,T}(\partial M)$, $t>0$
\begin{align}
\frac{|\nabla u|^2}{u^2} &+\delta\frac{H(u)}{u} -\delta \frac{u_t}{u}
 \le \frac{m\delta^2}{2(1-\epsilon)}\Big(\frac{1}{t}+ \frac{k}{\delta-1}+ \widetilde{\Phi} +\mu^H\Big) +  \frac{\delta}{p} \sqrt{ \frac{m}{2(1-\epsilon)}}\gamma^H
\end{align}
on $\partial M \times [T_0-T, T_0]$.  Here  $\epsilon \in (0,1)$, $\delta>1$ and  $p>0$ such that 
$$\frac{2(1-\epsilon)}{m\delta p}\ge \frac{1}{\epsilon}-1+\frac{(1-\delta)^2}{8}$$
and
$$\widetilde{\Phi} = \frac{1}{R^2}\Big[\Big(\frac{m\delta^2}{4(1-\epsilon)(\delta p+\delta-1)}\Big)C_1^2 +C_2 +  C_1 [(m-1)kR+l]R\Big]\Phi\Big]$$
for some constants positive $C_1$ and $C_2$.
Furthermore,
$$\mu^H : =   \sup_{\mathcal{Q}_{2R,T} } \{[H_u]_+\}  \ \ \text{and}\ \ \gamma^H := \sup_{\mathcal{Q}_{2R,T} }  \{[\Omega^H]_+\},$$ 
where $[H_u]_+ = [H'(u)-H(u)/u]_+$ and $[\Omega^H]_+ = [H'(u)-H(u)/u-(\alpha/p^2)uH''(u)]_+$.
\end{theorem}

In the above results the lower bound on the weighted Bakry-\'Emery Ricci tensor, $\mathscr{R}ic_\phi^m\ge -(m-1)k$, concerns the infimum of     
$\mathscr{R}ic_\phi^m$ only in the interior of $M$, whilst the lower bound on the weighted mean curvature of the boundary implies some weakconvexity condition on the boundary $\partial M$.     In \cite{[LY86]} Li and Yau have considered gradient estimates on linear heat equation with Neumann boundary condition on manifold whose boundary is convex in the sense that the second fundamental form of the boundary is nonnegative. (See also Wang \cite{[Wang97]}, Chen \cite{[Chen]} and Ramos Oliv\'e \cite{[Olive]}).  In some recent works Dung, Dung and Wu \cite{[DDW21]} studied Souplet-Zhang type gradient estimates on weighted heat equation with Dirichlet boundary condition on a smooth metric space with compact boundary satisfying $\mathscr{H}_\phi\ge -l$ based on the idea of Kunikawa and Sakurai \cite{[KS21]}. Fu and Wu \cite{[FW21]} extended this result to nonlinear parabolic equation which is somewhat a special case of ours with $H(u)=au\log u$, $a\in \mathbb{R}$.  In \cite{[Ab22]},  we established Souplet-Zhang type gradient estimates on a generalized nonlinear parabolic equation smooth metric measure spaces with compact boundary. A generalized Laplacian comparison result (Theorem \ref{thmW})  plays a key role in the proofs of all these results. The Sakurai generalized Laplacian comparison result is quite different from Wei-Wylie comparison \cite{[WeW09]} without boundary.  Fu and Wu \cite{[FW21]} remarked that Wei-Wylie comparison theorem requires the radius of geodesic ball $R\ge R_0$ for some constant $R_0>0$.

\section{Hamilton-Souplet-Zhang  type gradient estimates I}\label{sec3}
The proof of the following theorem is presented in this section.
\begin{theorem}\label{thm11}
Let $M=(M,g,e^{-\phi}dv)$ be an $n$-dimensional complete smooth metric measure space with compact boundary satisfying the bounds
\begin{align*}
\mathscr{R}ic_\phi^m \ge -(m-1)\kappa \ \ \text{and}\ \ \mathscr{H}_\phi \ge - \ell
\end{align*}
for some constants $\kappa, \ell \ge 0$. Let $u=u(x,t)$ be a smooth positive solution to \eqref{11} (i.e.   $(\partial_t -\mathscr{L}_\phi) u =H(u)$ in $\mathcal{Q}_{R,T}(\partial M$). Suppose further that $u$ satisfies the Dirichlet boundary condition (i.e., $u(\cdot,t)|_{\partial M}$ is constant for each time $t \in (-\infty,+\infty)$),
\begin{align*}
u_\nu\ge 0\ \ \text{and}\ \ u_t \le H(u)
\end{align*}
over $\partial M \times [T_0-T, T_0]$. Then there exists a constant $C>0$ depending only on $n,m$, such that for all $(x,t)\in \mathcal{Q}_{R/2,T/2}(\partial M)$
\begin{align}
\frac{|\nabla u|}{1+\sqrt{u}} & \le C\Bigg\{\left( \frac{1}{R} + \frac{1}{\sqrt{T}}+\sqrt{L}\right)\Big(\sup_{\mathcal{Q}_{R,T}(\partial M) }\{u\}\Big) \nonumber\\
& \hspace{2cm} + \Big(\sup_{\mathcal{Q}_{R,T}(\partial M) }\{\sqrt{\mathbb{M}_+[u]}\}\Big) \Big(\sup_{\mathcal{Q}_{R,T}(\partial M) }\{\sqrt{u}\}\Big) \Bigg\} ,
\end{align}
where  $L= (\kappa^2+\ell^4)^{1/2}$,
$$\mathbb{M}_+[u] = \left[\frac{u[(1/2+\sqrt{u})\widehat{H}(u) - (1+\sqrt{u})H'(u)]}{(1+\sqrt{u})^2} \right]_+ \ \ \text{and}\ \ \widehat{H}(u):=\frac{H(u)}{u}.$$  
\end{theorem}

To prove the above theorem we need some fundamental lemmas which we state and prove below.
\begin{lemma}\label{lem1T1}
Let $f(x,t)=u^{\beta}(x,t)$, $\beta \in (0,1)$, be a smooth function where $u(x,t)$ is a positive smooth solution to \eqref{11}. Then $f$ satisfies the nonlinear equation
\begin{align}\label{eqTL1}
(\mathscr{L}_\phi -\partial_t)f = \beta^{-1}(\beta-1)|\nabla f|^2/f -\beta f \widehat{H}(f),
\end{align}
where $ \widehat{H}(f) = H(f^{1/\beta})/f^{1/\beta} = H(u)/u.$
\end{lemma}

\proof
By direct computation for $f=u^{\beta}$ we have
$f_t=\beta u^{\beta-1}u_t$, $\nabla f=\beta u^{\beta-1}\nabla u$,  $|\nabla f|^2/f^2 = \beta^2|\nabla u|^2/u^2$ and
$\Delta f = \beta u^{\beta-1}\Delta u+ \beta(\beta-1)u^{\beta-2}|\nabla u|^2.$ Therefore
\begin{align*}
\mathscr{L}_\phi f & = \Delta f - \nabla \phi \nabla f\\
& = \beta u^{\beta-1}\mathscr{L}_\phi u +\beta(\beta-1)u^{\beta-2}|\nabla u|^2\\
& = \beta u^{\beta-1}(u_t-H(u)) +(\beta-1)u^{\beta-2}|\nabla u|^2\\
& = f_t - \beta f H(f^{1/\beta})/f^{1/\beta} + \beta^{-1}(\beta-1)|\nabla f|^2/f
\end{align*}
and consequently we arrive at \eqref{eqTL1}.

\qed

\begin{lemma}\label{lem2T1}
Let $f(x,t)=u^{\beta}(x,t)$, $\beta \in (0,1)$, be a smooth function where $u(x,t)$ is a positive smooth solution to \eqref{11}.  Set $w:=|\nabla f|^2/(1+f)^2$, then $w$ satisfies
\begin{align}\label{eqTL2}
(\mathscr{L}_\phi -\partial_t)w  = & \frac{|\nabla^2 f|^2}{(1+f)^2} + \frac{2\mathscr{R}ic_\phi(\nabla f, \nabla f)}{(1+f)^2}  - \frac{4\nabla f\nabla(|\nabla f|^2)}{(1+f)^3} + \frac{6|\nabla f|^4}{(1+f)^4} \nonumber\\
& -2\beta^{-1}(\beta-1) \left[\frac{|\nabla f|^4}{(1+f)^3f}  - \frac{\nabla f \nabla(|\nabla f|^2)}{(1+f)^2f} + \frac{|\nabla f|^4}{(1+f)^2f^2}\right] \nonumber\\
&  +2\beta\left[ \frac{|\nabla f|^2 f\widehat{H}(f)}{(1+f)^3} - \frac{\nabla f\nabla(f\widehat{H}(f))}{(1+f)^2} \right],
\end{align}
where $ \widehat{H}(f) = H(f^{1/\beta})/f^{1/\beta}.$
\end{lemma}

\proof
From the lemma hypothesis $w$ is given by $w=|\nabla f|^2/(1+f)^2$ with $f=u^\beta$, $0<\beta<1$. Based on this we can directly compute $\mathscr{L}_\phi w$ and $\partial_t w$. For this course, we first note that
\begin{align}\label{eqTL3}
\partial_t w = \partial_t\left(|\nabla f|^2/(1+f)^2\right) = \frac{2\nabla f\nabla f_t}{(1+f)^2} - \frac{2|\nabla f|^2f_t}{(1+f)^3}.
\end{align}
Similarly, we compute the gradient and Laplacian of $w$ as follows
\begin{align*}
\nabla w & = \frac{\nabla |\nabla f|^2}{(1+f)^2} - \frac{2|\nabla f|^2\nabla f}{(1+f)^3},\\
\Delta w &= \frac{\Delta |\nabla f|^2}{(1+f)^2} - \frac{4 \nabla|\nabla f|^2\nabla f}{(1+f)^3} - \frac{2|\nabla f|^2\Delta f}{(1+f)^3}  + \frac{6|\nabla f|^4}{(1+f)^4}
\end{align*}
which then yields
\begin{align}\label{eqTL4}
\mathscr{L}_\phi w & = \Delta w - \nabla \phi\nabla w \nonumber\\
& =  \frac{\mathscr{L}_\phi(|\nabla f|^2)}{(1+f)^2} - \frac{4 \nabla|\nabla f|^2\nabla f}{(1+f)^3} - \frac{2|\nabla f|^2 \mathscr{L}_\phi f}{(1+f)^3}  + \frac{6|\nabla f|^4}{(1+f)^4}.
\end{align}
Now putting together the individual descriptions of $\partial_t w$ and $\mathscr{L}_\phi w$ from \eqref{eqTL3} and \eqref{eqTL4} respectively, and taking into account the relevant cancellations, we get 
\begin{align*}
(\mathscr{L}_\phi - \partial) w & =  \frac{\mathscr{L}_\phi(|\nabla f|^2)}{(1+f)^2} - \frac{4 \nabla|\nabla f|^2\nabla f}{(1+f)^3} - \frac{2|\nabla f|^2 (\mathscr{L}_\phi f-f_t)}{(1+f)^3}  + \frac{6|\nabla f|^4}{(1+f)^4} - \frac{2\nabla f\nabla f_t}{(1+f)^2}.
\end{align*}
Combining the last equation with the weighted Bochner formula \eqref{e1}, we arrive at 
\begin{align}\label{eqTL6}
(\mathscr{L}_\phi - \partial) w  =  &  \frac{2|\nabla^2 f|^2}{(1+f)^2} + \frac{2\mathscr{R}ic_\phi(\nabla f, \nabla f)}{(1+f)^2}  +  \frac{2\nabla f (\mathscr{L}_\phi f -\nabla f_t)}{(1+f)^2} \nonumber  \\
&  - \frac{2|\nabla f|^2 (\mathscr{L}_\phi f-f_t)}{(1+f)^3} - \frac{4 \nabla|\nabla f|^2\nabla f}{(1+f)^3}+ \frac{6|\nabla f|^4}{(1+f)^4}.
\end{align}
Referring to \eqref{eqTL1} we compute directly
\begin{align}\label{eqTL7}
\nabla f\nabla (\mathscr{L}_\phi f - f_t) &= \nabla f\nabla (\beta^{-1}(\beta-1)|\nabla f|^2/f -\beta f \widehat{H}(f))\nonumber\\
& =  \beta^{-1}(\beta-1)\left[\nabla f \nabla(|\nabla f|^2)/f - |\nabla f|^2/f^2\right] - \beta \nabla f \nabla (f\widehat{H}(f)).
\end{align}
Finally, plugging in  \eqref{eqTL1} and  \eqref{eqTL7} into  \eqref{eqTL6}, and rearranging terms gives the desired conclusion.

\qed

\begin{lemma}\label{lem3T1}
Under the assumptions of Lemma \ref{lem2T1}, suppose $\mathscr{R}ic_\phi^m \ge -(m-1)\kappa$, $\kappa\ge 0$, then 
\begin{align}\label{eqTL8}
(\mathscr{L}_\phi -\partial_t)w & \ge -\frac{2[f-\beta^{-1}(\beta-1)(1+f)]}{1+f} \langle \nabla f, \nabla w\rangle\nonumber \\
& \hspace{1cm}  -2\beta^{-1}(\beta-1)\frac{(1+f)}{f} w^2 - 2(m-1)\kappa w + 2\mathbb{M}_\beta(f)w,
\end{align}
where 
$$ \mathbb{M}_\beta(f):= \left[\frac{(1-\beta+f)\widehat{H}(f)-(1+f)H'(f)}{1+f}\right] \ \  \text{and} \ \ \widehat{H}(f) = H(f^{1/\beta})/f^{1/\beta}.$$
\end{lemma}

\proof 
The proof is based on further simplification of \eqref{eqTL2} and application of weighted  $m$-Ricci tensor lower bound.  First recall that $\widehat{H}(f) = H(f^{1/\beta})/f^{1/\beta}$, so that 
\begin{align*}
 \nabla (f\widehat{H}(f)) & = \nabla ( H(f^{1/\beta}) f^{-1/\beta+1})\\
 &  =  (1/\beta)\nabla f H'(f^{1/\beta}) - (1/\beta - 1) \nabla f H(f^{1/\beta})/ f^{1/\beta}
\end{align*}
and then
\begin{align*}
2\beta \frac{\nabla f \nabla (f\widehat{H}(f))}{(1+f)^2} & = 2 \frac{|\nabla f|^2}{(1+f)^2} H'(f^{1/\beta}) + 2(\beta-1)  \frac{|\nabla f|^2}{(1+f)^2} \widehat{H}(f^{1/\beta})\\
& = 2 \left[H'(f^{1/\beta}) + (\beta-1) \widehat{H}(f^{1/\beta}) \right] w.
\end{align*}
By this, it can be seen clearly that 
\begin{align}\label{eqTL9}
2\beta\left[ \frac{|\nabla f|^2 f\widehat{H}(f)}{(1+f)^3} - \frac{\nabla f\nabla(f\widehat{H}(f))}{(1+f)^2} \right] & = 2 \left[\Big(\frac{\beta f}{(1+f)}-(\beta-1)\Big)\widehat{H}(f^{1/\beta}) - H'(f^{1/\beta})\right]w\nonumber\\
& = 2 \mathbb{M}_\beta(f) w,
\end{align}
where  $\displaystyle \mathbb{M}_\beta(f):= \left[\frac{(1-\beta+f)\widehat{H}(f^{1/\beta})-(1+f)H'(f^{1/\beta})}{1+f}\right].$

\noindent
Next considering the following identity (which can be verified by direct derivation)
\begin{align*}
\frac{4 \nabla^2 f |\nabla f|^2 }{(1+f)^3} = \frac{2 \nabla f \nabla(|\nabla f|^2)}{(1+f)^3}
\end{align*}
one can then simplify some other terms of \eqref{eqTL2}  as follows
\begin{align}\label{eqTL10}
\frac{2|\nabla^2 f|^2}{(1+f)^2}- &\frac{4 \nabla|\nabla f|^2\nabla f}{(1+f)^3}+ \frac{6|\nabla f|^4}{(1+f)^4}\nonumber\\
& =  \frac{2}{(1+f)^2}\left| \nabla^2f - \frac{\nabla f\otimes \nabla f}{1+f}\right|^2 - \frac{2 \nabla|\nabla f|^2\nabla f}{(1+f)^3}+ \frac{4|\nabla f|^4}{(1+f)^4}.
\end{align}
Refering to the definition of  weighted $m$-Ricci tensor given in the introduction  we write
\begin{align}\label{eqTL11}
\frac{2\mathscr{R}ic_\phi(\nabla f, \nabla f)}{(1+f)^2} &= \frac{2[(m-n)\mathscr{R}ic_\phi^m(\nabla f, \nabla f) + \nabla \phi \otimes \nabla \phi(\nabla f, \nabla f)}{(m-n)(1+f)^2} \nonumber\\
& \ge -2(m-1)\kappa w
\end{align}
since $\mathscr{R}ic_\phi^m(\nabla f, \nabla f) \ge - (m-1)\kappa|\nabla f|^2$ from the hypothesis and the observation that
$$ \frac{\nabla \phi \otimes \nabla \phi(\nabla f, \nabla f)}{(m-n)(1+f)^2} =  \frac{|\langle \nabla \phi, \nabla f\rangle|^2}{(m-n)(1+f)^2} \ge 0$$
for $m\ge n$ ($m=n$ only for $f$ constant). Putting \eqref{eqTL9}-\eqref{eqTL11} into \eqref{eqTL2} we have 
\begin{align}\label{eqTL12}
(\mathscr{L}_\phi -\partial_t)w  \ge & - \frac{2\nabla f\nabla(|\nabla f|^2)}{(1+f)^3} + \frac{4|\nabla f|^4}{(1+f)^4}  -2\beta^{-1}(\beta-1) \Bigg[\frac{|\nabla f|^4}{(1+f)^3f} \nonumber\\
&  - \frac{\nabla f \nabla(|\nabla f|^2)}{(1+f)^2f} + \frac{|\nabla f|^4}{(1+f)^2f^2}\Bigg] -2(m-1)\kappa w  + 2\mathbb{M}_\beta(f)w. 
\end{align}
To simplify further some terms of \eqref{eqTL12} we take into consideration the following identity
$$2\nabla f\nabla w = \frac{2\nabla f\nabla(|\nabla f|^2)}{(1+f)^2} - \frac{4|\nabla f|^4}{(1+f)^3}$$
so that we compute the following identities
\begin{align}\label{eqTL13}
 - \frac{2\nabla f\nabla(|\nabla f|^2)}{(1+f)^3} + \frac{4|\nabla f|^4}{(1+f)^4} = - \frac{2}{1+f}\nabla f\nabla w
\end{align}
and 
\begin{align}\label{eqTL14}
2\beta^{-1}(\beta-1) \left[ \frac{\nabla f \nabla(|\nabla f|^2)}{(1+f)^2f} - \frac{|\nabla f|^4}{(1+f)^3f}\right] =  2\beta^{-1}(\beta-1)\frac{1}{f} \nabla f\nabla w.
\end{align}
To this end, \eqref{eqTL13} and \eqref{eqTL14} show that the first four terms on the right hand side of  \eqref{eqTL12} are simplified to
\begin{align}\label{eqTL15}
- \frac{2}{1+f}\nabla f\nabla w + 2\beta^{-1}(\beta-1)\frac{1}{f} \nabla f\nabla w = -\frac{2[f-\beta^{-1}(\beta-1)(1+f)]}{1+f}  \nabla f\nabla w.
\end{align}
Replacing the first four terms on the right hand side of  \eqref{eqTL12} by \eqref{eqTL15} gives the required inequality which is \eqref{eqTL8}.

\qed

\begin{remark}\label{rem1T1}
Informed by the sign of the constant $\beta^{-1}(\beta-1)$ appearing in the inequality \eqref{eqTL8} which is always negative in the range $0<\beta<1$, we can normalize the constant as $\beta^{-1}(\beta-1)=1$. This gives a good impression that $\beta$ can be chosen to be $\beta=1/2$ in the later analysis.
\end{remark}

\subsection{proof of Theorem \ref{thm11}}

\proof
Suppose $f=u^\beta$, $\beta \in (0,1)$ and $w= |\nabla \log(1+f)|^2$ satisfying the assumptions of Lemma \ref{lem1T1}. In the case $w$ is constant on $\mathcal{Q}_{R/2,T}(\partial M)$ there remain nothing to prove. To this end we assume that $w$ is non-constant in $\mathcal{Q}_{R/2,T}(\partial M)$. Let there exists a smooth cut-off function $\psi$ verifying those properties listed in Lemma \ref{lem23}. Clearly $\psi w$ is positive and attains its maximum at a point in $\mathcal{Q}_{R/2,T}(\partial M)$, say $(x_1,t_1)$. Our argument will henceforth split into two cases, according to whether $x_1 \not\in \partial M$ or $x_1 \in \partial M$. 

\noindent
{\bf Case 1:} If $x_1 \not\in \partial M$, we may assume without loss of generality that $x_1 \not\in \partial M\cup \ \text{cut}(\partial M)$ (where $\text{cut}(\partial M)$ is the cut locus of $\partial M$) by invoking the well known's Calabi's trick (see \cite{[Cal],[LY86],[SY94]}). Then at $(x_1,t_1)$ we have 
$$\nabla(\psi w)=0, \ \  \mathscr{L}_\phi(\psi w) \le 0 \ \ \text{and} \ \   \partial_t(\psi w) \ge 0.$$
To start with we  first observe that an elementary calculation yields
\begin{align}\label{A1}
(\mathscr{L}_\phi - \partial_t)(\psi w) & = \psi (\mathscr{L}_\phi- \partial_t)w + 2 \nabla \psi \nabla w + w(\mathscr{L}_\phi - \partial_t)\psi.
\end{align}
Combining Lemma \ref{lem1T1} with  \eqref{A1} while denoting
$A:=2[(f-\beta^{-1}(\beta-1)(1+f))/(f(1+f))]\nabla f$ and 
$B:=2\beta^{-1}(\beta-1)(1+f)/f^2$ for the sake of convenience and brevity, we obtain
\begin{align}\label{A2}
(\mathscr{L}_\phi - \partial_t)(\psi w)  &\ge  -A\psi\nabla w - B\psi w^2 +2\mathbb{M}_\beta[f]\psi w - 2(m-1)k\psi w \nonumber \\
& \hspace{1cm} + 2 \nabla \psi \nabla w + w(\mathscr{L}_\phi - \partial_t)\psi.
\end{align}
Proceeding forward, we note the following basic vector identities using Leibniz rule: 
$$A\psi\nabla w = A\nabla(\psi w)-(A\nabla \psi)w \ \ \text{and}\ \ \nabla\psi\nabla w = (\nabla \psi/\psi)\nabla(\psi w) -(|\nabla\psi|^2/\psi)w.$$
Thus upon substitution of these identities into \eqref{A2} we have 
\begin{align}\label{A3}
(\mathscr{L}_\phi - \partial_t)(\psi w)  \ge  & -A\nabla(\psi w) +(A\nabla \psi)w - B\psi w^2  +2\frac{\nabla \psi}{\psi} \nabla(\psi w) - 2\frac{|\nabla\psi|^2}{\psi}w   \nonumber \\
& +2\mathbb{M}_\beta[f]\psi w + 2 \nabla \psi \nabla w + w[\mathscr{L}_\phi - \partial_t -  2(m-1)k]\psi.
\end{align}
Rearranging \eqref{A3} and making use of the fact that localized function $\psi w$ attains its maximum at $(x_1,t_1)$,  where we note that $\nabla(\psi w)=0$ and $(\mathscr{L}_\phi - \partial_t)(\psi w)\le 0$,  and obtain
\begin{align}\label{A4}
-B\psi w^2 \le (A\nabla \psi)w + 2\frac{|\nabla\psi|^2}{\psi}w - 2\mathbb{M}_\beta[f]\psi w -  w[\mathscr{L}_\phi - \partial_t -  2(m-1)k]\psi.
\end{align}
Since $B:=2\beta^{-1}(\beta-1)(1+f)/f^2$ and $0<\beta<1$, we note that $\beta^{-1}(\beta-1)<0$, meaning that $B$ is negative and $-B>0$. So we multiply through \eqref{A4} by $(-B)^{-1}$, and we then have
\begin{align}\label{A5}
\psi w^2& \le \left(\frac{A}{B}\nabla \psi\right)w - \beta(\beta-1)^{-1}\frac{|\nabla\psi|^2}{\psi}\frac{f^2}{1+f}w + \beta(\beta-1)^{-1}\mathbb{M}_\beta[f]\frac{f^2}{1+f}\psi w\nonumber \\
& \hspace{1cm} + \beta(\beta-1)^{-1}\frac{|\nabla\psi|^2}{\psi}\frac{f^2}{2(1+f)}w (\mathscr{L}_\phi - \partial_t -  2(m-1)k)\psi.
\end{align}
The next step is to apply inequality \eqref{A5} to derive estimates at $(x,t)$. Here repeated applications of Cauchy-Schwarz inequality, Young'sinequality and the properties of the cut-off function $\psi$ presented in Lemma \ref{lem23} will give appropriate bounds for each term on the right hand side of \eqref{A5}. Before we proceed, note again that $u^\beta=f>0$, $0<\beta<1$ implies
$0<f/(1+f)\le 1$, \  $\beta^{-1}(\beta-1)<0,$ \ $\beta(\beta-1)^{-1}<0 $ and $f/(1+f) -\beta(\beta-1)^{-1}>0$.

\noindent
By definition
$$\frac{A}{B} = \beta(\beta-1)^{-1}\left[\frac{f - \beta^{-1}(\beta-1)(1+f)}{1+f}\right] \frac{f\nabla f}{(1+f)}.$$
Since the quantity in front  of $(f\nabla f)/(1+f)$ in the last equality is negative we estimate the first term on the right hand side of \eqref{A5} as follows:
\begin{align}\label{A6}
 \left(\frac{A}{B}\nabla \psi\right)w  & = \beta(\beta-1)^{-1}\left[\frac{f - \beta^{-1}(\beta-1)(1+f)}{1+f}\right] \frac{f\nabla f}{(1+f)}\nabla \psi w \nonumber \\
 & \le  \beta(1-\beta)^{-1}\left[\frac{f - \beta^{-1}(\beta-1)(1+f)}{1+f}\right]f \frac{|\nabla f|}{(1+f)}| \nabla \psi| w \nonumber \\
 & \le  (1-\beta)^{-1} \Big(\sup_{\mathcal{Q}_{R,T}} f\Big)|\nabla \psi| w^{3/2} \nonumber \\
&  \le \frac{1}{8}\psi w^2 + C_\beta \left(\frac{|\nabla \psi| \Big(\sup_{\mathcal{Q}_{R,T}} f\Big)}{\psi^{3/4}} \right)^4 \nonumber \\
& \le \frac{1}{8}\psi w^2 + \frac{C_\beta}{R^4} \Big(\sup_{\mathcal{Q}_{R,T}} f\Big)^4,
\end{align}
where $\mathcal{Q}_{R,T}:=\mathcal{Q}_{R,T}(\partial M)$ and $C_\beta$ is a constant which depends only on $\beta$,  $n$ and $m$, the value of which may vary from line to line as remarked before. (Similarly any $C_\beta$ appearing in the rest analysis depends only on $\beta$,  $n$ and $m$ with varying values).  

Estimating the second term on the right hand side of \eqref{A5} we proceed in a likewise manner as follows
\begin{align}\label{A7}
- \beta(\beta-1)^{-1}\frac{|\nabla\psi|^2}{\psi}\frac{f^2}{1+f}w & =   \beta(1-\beta)^{-1}\frac{|\nabla\psi|^2}{\psi}\frac{f^2}{1+f}w  \nonumber\\
& \le \beta(1-\beta)^{-1}\psi^{1/2} \frac{|\nabla\psi|^2}{\psi^{3/2}}\Big(\sup_{\mathcal{Q}_{R,T}} f\Big)^2\nonumber\\
& \le \frac{1}{8}\psi w^2 + C_\beta\Big( \frac{|\nabla\psi|^2}{\psi^{3/2}}\Big)^2\Big(\sup_{\mathcal{Q}_{R,T}} f\Big)^4\nonumber\\
& \le \frac{1}{8}\psi w^2 + \frac{C_\beta}{R^4} \Big(\sup_{\mathcal{Q}_{R,T}} f\Big)^4.
\end{align}
For the third term on the right hand side of \eqref{A5} we  first write $\bar{\mathbb{M}}_{\beta} [f]= f^2\mathbb{M}_\beta [f]/(1+f)$, ( where $\mathbb{M}_\beta[f]$ is as defined in Lemma \ref{lem3T1})  and then estimate as follows:
\begin{align}\label{A8}
\beta(\beta-1)^{-1}\mathbb{M}_\beta[f]\frac{f^2}{1+f}\psi w & =  \beta(\beta-1)^{-1}\bar{\mathbb{M}}_\beta[f]\psi w \nonumber\\
& \le  \beta(1-\beta)^{-1}\bar{\mathbb{M}}_{\beta+}[f]\psi w \nonumber\\
& \le \frac{1}{8}\psi w^2 + C_\beta \Big(\sup_{\mathcal{Q}_{R,T}} \bar{\mathbb{M}}_{\beta+}[f]\Big)^2,
\end{align}
where 
$$\bar{\mathbb{M}}_{\beta+}[f] := \left[\frac{f^2[(1-\beta+f)\widehat{H}(f^\frac{1}{\beta})  -(1+f)H'(f^\frac{1}{\beta})]}{(1+f)^2} \right]_+.
$$
We break the fourth term to individual item and estimate each of them and later sum up their estimates.
For the term with time derivative of $\psi$:
\begin{align}\label{A9}
-\beta(\beta-1)^{-1}\frac{f^2}{2(1+f)}w\psi_t &= \beta (1-\beta)^{-1}\frac{f^2}{2(1+f)}w|\psi_t|\nonumber\\
& \le \frac{1}{2}\beta (1-\beta)^{-1}\psi^{\frac{1}{2}} w \frac{|\psi_t|}{\psi^{\frac{1}{2}}}\Big(\sup_{\mathcal{Q}_{R,T}} f\Big)^2\nonumber\\
& \le \frac{1}{8}\psi w^2 + C_\beta \frac{|\psi_t|^2}{\psi} \Big(\sup_{\mathcal{Q}_{R,T}} f\Big)^4\nonumber\\
& \le \frac{1}{8}\psi w^2 + \frac{C_\beta}{T^2} \Big(\sup_{\mathcal{Q}_{R,T}} f\Big)^4.
\end{align}
Secondly, the term with $(m-1)\kappa$ is estimated in a similar way:
\begin{align}\label{A10}
-\beta(\beta-1)^{-1}(m-1)\kappa\psi w\frac{f^2}{1+f} &= \beta (1-\beta)^{-1} (m-1)\kappa\psi w\frac{f^2}{1+f}  \nonumber\\
&\le   \beta (1-\beta)^{-1} (m-1)\kappa  (\psi^{\frac{1}{2}}w) \psi^{\frac{1}{2}}\Big(\sup_{\mathcal{Q}_{R,T}} f\Big)^4  \nonumber\\
& \le \frac{1}{8}\psi w^2 + C_\beta\kappa^2 \Big(\sup_{\mathcal{Q}_{R,T}} f\Big)^4.
\end{align}
Lastly, we compute the term involving the weighted Laplacian of $\psi$ as follows. Here we apply the weighted Laplacian comparison $\mathscr{L}_\phi \varrho_{\partial M}(x) \le (n-1)\kappa R + \ell$, $\kappa,\ell \ge 0$, for all $x \in B_R(\partial M)$ since $\mathscr{R}ic_\phi^m \ge -(m-1)\kappa$ and $\mathscr{H}_\phi \ge -\ell$ (that is, the weighted Ricci tensor and the boundary mean curvature bounds (see Theorem \ref{thmW})):
\begin{align}\label{A11}
\frac{\beta f^2}{2(\beta-1)(1+f)}w\mathscr{L}_\phi\psi &=-  \frac{\beta f^2}{2(1-\beta)(1+f)} \left[ \psi_\varrho \mathscr{L}_\phi\varrho + \psi_{\varrho \varrho} |\nabla \varrho|^2 \right]w  \nonumber\\
&\le \frac{\beta f^2}{2(1-\beta)(1+f)} \left[| \psi_\varrho| ((n-1)\kappa R+\ell) + \psi_{\varrho \varrho} |\nabla \varrho|^2 \right]w  \nonumber\\ 
&\le \frac{\beta}{2(1-\beta)} \psi^{1/2}w\left[ \frac{|\psi_{\varrho \varrho}|}{\psi^{1/2}}     + \frac{|\psi_\varrho|}{\psi^{1/2}}\Big((n-1)\kappa R+ \ell\Big) \right]\Big(\sup_{\mathcal{Q}_{R,T}} f\Big)^2 \nonumber\\
&\le \frac{1}{8}\psi w^2 + C_\beta \left[ \frac{|\psi_{\varrho \varrho}|^2}{\psi}     + \frac{|\psi_\varrho|^2}{\psi} (n^2\kappa^2R^2 +l^2) \right]\Big(\sup_{\mathcal{Q}_{R,T}} f\Big)^4\nonumber\\
&\le \frac{1}{8}\psi w^2 + C_\beta \left[ \frac{1}{R^4}  +\kappa^2 +\ell^4\right]\Big(\sup_{\mathcal{Q}_{R,T}} f\Big)^4.
\end{align}
With bounds \eqref{A6}--\eqref{A11} at our disposal, we continue the analysis for $x\not\in\partial M \cup\ \text{cut}(\partial M)$ by substituting the bounds back into \eqref{A5}, simplifying further and rearranging terms to get the following inequality
\begin{align*}
\psi w^2& \le\frac{3}{4}\psi w + C_\beta\left\{ \Big(\frac{1}{R^4}+\frac{1}{T^2}+L^2\Big)\Big(\sup_{\mathcal{Q}_{R,T}} f\Big)^4 + \Big(\sup_{\mathcal{Q}_{R,T}} \mathbb{M}_{\beta+}[f]\Big)^2 \right\},
\end{align*}
where $L=(\kappa^2+\ell^4)^{1/2}$ at $(x_1,t_1)$.  
It therefore follows by reverting to the maximal characterization of the point $(x_1,t_1)$ that we have for all $(x,t)\in \mathcal{Q}_{R/2,T/2}(\partial M)$
\begin{align*}
\psi^2w^2(x,t) & \le \psi^2w^2(x_1,t_1)\le \psi w^2(x_1,t_1)\\
& \le C_\beta \left\{ \frac{1}{R^4}\Big(1+\frac{R^4}{T^2}+L^2R^4\Big)\Big(\sup_{\mathcal{Q}_{R,T}} f\Big)^4 + \Big(\sup_{\mathcal{Q}_{R,T}} \mathbb{M}_{\beta+}[f]\Big)^2 \right\}.
\end{align*}
Since $\psi(x,t)\equiv 1$ on $\mathcal{Q}_{R/2,T/2}(\partial M)$ and by definition $ w =|\nabla f|^2/(1+f)^2$ we obtain
\begin{align*}
\frac{|\nabla f|}{1+f} \le C_\beta \left\{ \frac{1}{R}\Big(1+\frac{R}{\sqrt{T}}+\sqrt{L}R\Big)\Big(\sup_{\mathcal{Q}_{R,T}} f\Big) + \Big(\sup_{\mathcal{Q}_{R,T}} \mathbb{M}_{\beta+}[f]\Big)^{1/2} \right\}.
\end{align*}
Considering the range of $\beta\in (0,1)$ and referring to the comments made in Remark \ref{rem1T1} above, we can choose $\beta=1/2$. Noting that $f=u^\beta$ with $\beta=1/2$ yields $f=\sqrt{u}$ and $|\nabla f| = |\nabla u|/2\sqrt{u}$, we finally arrive at 
\begin{align*}
\frac{|\nabla \sqrt{u}|}{1+\sqrt{u}} \le C \left\{ \frac{1}{R}\Big(1+\frac{R}{\sqrt{T}}+\sqrt{L}R\Big)\Big(\sup_{\mathcal{Q}_{R,T}} u\Big) + \Big(\sup_{\mathcal{Q}_{R,T}} \sqrt{\mathbb{M}_{+}[u]}\Big)\Big(\sup_{\mathcal{Q}_{R,T}} \sqrt{u}\Big)  \right\},
\end{align*}
which says the desired gradient estimate \eqref{e1} holds in the interior for all $(x,t)\in \mathcal{Q}_{R/2,T/2}(\partial M)$.
It remains to prove that the gradient estimate is still valid on the boundary. We therefore consider the case $x_1\in \partial M$.

\noindent
{\bf Case 2:} Assume $x_1\in \partial M$.  Indeed, at the maximum point  $(x_1,t_1)$ of $(\psi w)$, we have $(\psi w)_\nu\ge 0$ and consequently, 
$\psi_\nu w+\psi w_\nu \ge 0$, which 
in particular  gives
$$w_\nu\ge 0.$$
Since $w=|\nabla\log(1+f)|^2$,  and $\log(1+f)$ satisfies the Dirichlet boundary condition, we use the Reilly-type formula \eqref{Rei} and get
\begin{align}\label{A12}
0\le w_\nu &= (|\nabla\log(1+f)|^2)_\nu \nonumber\\
&= 2(\log(1+f))_\nu\left[\mathscr{L}_\phi(\log(1+f))-(\log(1+f))_\nu \mathscr{H}_\phi \right].
\end{align}
Notice also  that $f$  satisfies the Dirichlet boundary condition since $u$ satisfies the Dirichlet boundary condition,   and then $|\nabla f| =f_v= \beta u^\beta u_\nu/u \ge 0$ since $u>0$ and $u_v\ge 0$ on $\partial M$.  By smoothness of the function $f$ we have
\begin{align}\label{A13}
\mathscr{L}_\phi(\log(1+f)) &= \Delta (\log(1+f)) -\nabla\phi  \nabla(\log(1+f))\nonumber\\
& = -\frac{|\nabla f|^2}{(1+f)^2} + \frac{\Delta f}{1+f} -  \frac{1}{1+f} \nabla\phi\nabla f \nonumber \\
& = -\frac{|\nabla f|^2}{(1+f)^2} + \frac{\mathscr{L}_\phi f}{1+f}.
\end{align}
Setting $\beta=1/2$, we have from Lemma \ref{lem1T1} that 
\begin{align}\label{A14}
f_t = \mathscr{L}_\phi f + |\nabla f|^2/f + f\widehat{H}(f)/2.
\end{align}
Thus 
\begin{align}\label{A15}
\mathscr{L}_\phi(\log(1+f)) =- \frac{(1+2f)}{f}(\log(1+f))_\nu^2 +  \frac{f_t -\frac{1}{2}f\widehat{H}(f)}{1+f}
\end{align}
by substituting \eqref{A14} into \eqref{A13} and then simplifying the resulting expression.
Noting also that on $\partial M$, 
$(\log(1+f))_\nu = f_v/(1+f) = |\nabla f|/(1+f) = \sqrt{w}$. 
It then follows from \eqref{A2} that 
\begin{align*}
0\le 2\sqrt{w}\left(-\frac{(1+2f)}{f}w +  \frac{f_t -\frac{1}{2}f\widehat{H}(f)}{1+f} - \sqrt{w}\mathscr{H}_\phi \right)
\end{align*}
which upon rearranging and considering the assumption  $u_t\le H(u)$ over $\partial M\times [T_0-T,T_0]$ yields $f_t-\frac{1}{2}f\widehat{H}(f)\le 0$
 and the inequality
 \begin{align}\label{A16}
\frac{(1+2f)}{f}w^{3/2} + \mathscr{H}_\phi w \le 0
 \end{align}
at $(x_1,t_1)$. Taking into account that $(1+2f)/f>1/f>0$, \eqref{A16} is then solved at $(x_1,t_1)$ to give
$$w(x_1,t_1)=0 \ \ \ \text{or} \ \ \ w^{1/2}(x_1,t_1)\le \ell \cdot \sup_{\mathcal{Q}_{R,T}} f$$
since $\mathscr{H}_\phi \ge -\ell$ by the theorem hypothesis. This thus means
$$\psi w(x_1,t_1)=0 \ \ \ \text{or} \ \ \ \psi w(x_1,t_1)\le \ell^2 \cdot (\sup_{\mathcal{Q}_{R,T}} f)^2$$
on $\mathcal{Q}_{R,T}(\partial M)$. 
Note that $\psi w=0$ indicates that $f$ (respectively $u$) is a constant from where the conclusion follows at once, whereas the case $(\psi w) \le \ell^2(\sup f)^2$ for all $(x,t) \in \mathcal{Q}_{R/2,T/2}(\partial M)$ together with the property $\psi(x,t)\equiv 1$ implies
\begin{align*}
\frac{|\nabla f|^2}{(1+f)^2}(x,t) = w(x,t)=\psi w(x,t)\le (\psi W)(x_1,t_1)\le \ell^2 \cdot (\sup_{\mathcal{Q}_{R,T}} f)^2.
\end{align*}
Taking into consideration that $f=\sqrt{u}$ and $|\nabla f|/(1+f) = |\nabla \sqrt{u}|/(1+\sqrt{u})$,  the last inequality implies the conclusion.

\qed

\section{Hamilton-Souplet-Zhang type gradient estimates II}\label{sec4}
In this section we discuss the proof of Theorem \ref{thm22} by stating and proving the following theorem.
\begin{theorem}\label{thm21}
Let $M=(M,g,e^{-\phi}dv)$ be an $n$-dimensional complete smooth metric measure space with compact boundary satisfying the bounds
\begin{align*}
\mathscr{R}ic_\phi^m \ge -(m-1)\kappa \ \ \text{and}\ \ \mathscr{H}_\phi \ge - \ell
\end{align*}
for some constants $\kappa, \ell \ge 0$. Let $u=u(x,t)$ be a smooth positive solution to \eqref{11}  in $\mathcal{Q}_{R,T}(\partial M)$. Suppose further that $u$ satisfies the Dirichlet boundary condition (i.e., $u(\cdot,t)|_{\partial M}$ is constant for each time $t \in (-\infty,+\infty)$),
\begin{align*}
u_\nu\ge 0\ \ \text{and}\ \ u_t \le H(u)
\end{align*}
over $\partial M \times [T_0-T, T_0]$. Then for $\varepsilon \in (0,  1/2)$, there exists a constant $C_\varepsilon>0$ depending only on $n,m, \varepsilon$, such that for all $(x,t)\in \mathcal{Q}_{R/2,T/2}(\partial M)$
\begin{align}\label{ethm21}
\frac{|\nabla u|}{u^{1-3\varepsilon/2}} & \le C_\varepsilon\Bigg\{\frac{1}{R} + \frac{1}{\sqrt{T}}+\sqrt{L}+ \Gamma^\varepsilon_+[u] \Bigg\}\Big(\sup_{\mathcal{Q}_{R,T}(\partial M) }\{\sqrt{u^{3\varepsilon}}\}\Big) ,
\end{align}
where  $L= (\kappa^2+\ell^4)^{1/2}$,
$$\Gamma^\varepsilon_+[u] =  \sup_{\mathcal{Q}_{R,T}(\partial M) }\left\{\left[ 2H'(u)+(3\varepsilon-2)\widehat{H}(u) \right]_+ \right\}^{1/2} \ \ \text{and}\ \ \widehat{H}(u):=\frac{H(u)}{u}.$$  
\end{theorem}

\begin{remark}
The proof of Theorem \ref{thm22} follows from the above theorem by choosing $\varepsilon = 1/3$ in the range $0<\varepsilon <1/2$.

On the other hand, passing to the limit $\varepsilon \searrow 0$ in \eqref{ethm21} we have the following estimates for all $(x,t)\in \mathcal{Q}_{R/2,T/2}(\partial M)$ under the assumption of Theorem \ref{thm21}:
\begin{align}
\frac{|\nabla u|}{u} & \le C_0\Bigg\{\frac{1}{R} + \frac{1}{\sqrt{T}}+\sqrt{L}+ \Gamma^0_+[u] \Bigg\}, 
\end{align}
where the constant $C_0>0$ depends only on $n$ and $m$,  and 
$$\Gamma^0_+[u] =  \sup_{\mathcal{Q}_{R,T}(\partial M) }\left\{\left[ 2(H'(u)-\widehat{H}(u)) \right]_+ \right\}^{1/2}.$$  
\end{remark}

Next we prove some fundamental lemmas which will be applied in the proof of Theorem \ref{thm21}.
\section*{Fundamental lemmas}
In this section,we define a qunatity $w$ as $w(x,t)= f(x,t)|\nabla f(x,t)|^2$ with $f=u^\varepsilon$ and $\varepsilon\in (0,1)$, where $u>0$ solves the nonlinear heat-type equation \eqref{11}.  Going by Lemma \ref{lem1T1} of the previous section, it is clear that $f_t=\varepsilon u^{\varepsilon-1}u_t$, $\nabla f=\varepsilon u^{\varepsilon-1}\nabla u$,  $|\nabla f|^2/f^2 = \varepsilon^2|\nabla u|^2/u^2$ and $f$ satisfies 
\begin{align}\label{e3.1}
f_t =  \mathscr{L}_\phi f - \varepsilon^{-1}(\varepsilon-1)|\nabla f|^2/f +\varepsilon f H(f^{1/\varepsilon})/f^{1/\varepsilon}
\end{align}
Furthermore, we will denote $\widehat{H}(f)=H(f^{1/\varepsilon})/f^{1/\varepsilon}$ and note that the quantity $\widehat{H}(f)$ is the same as $\widehat{H}(u)=H(u)/u$.  

\begin{lemma}\label{lem3.1}
Let $f=u^\varepsilon$,  $\varepsilon \in (0,1)$, be a smooth function where $u(x,t)$ is a positive smooth solution to \eqref{11}. Then $w =f|\nabla f|^2$ satisfies the nonlinear equation
\begin{align}\label{e3.2}
(\mathscr{L}_\phi -\partial_t)w  =&   2f|\nabla^2f|^2 + 2\varepsilon^{-1}(\varepsilon-1)\nabla f\nabla(|\nabla f|^2) - \varepsilon^{-1}(\varepsilon-1)|\nabla f|^4/f\nonumber\\
& + 2f\mathscr{R}ic_\phi(\nabla f,\nabla f) + 2\nabla f\nabla (|\nabla f|^2)-\varepsilon|\nabla f|^2f\widehat{H}(f)\\
&  - 2\epsilon f \nabla f \nabla (f\widehat{H}(f)).\nonumber
\end{align}
\end{lemma}

\proof
By direct computation 
\begin{align}\label{e3.3}
w_t=f_t |\nabla f|^2 + 2f \langle\nabla f,\nabla f_t\rangle.
\end{align}
Using the weighted Bochner formula 
\begin{align}\label{e3.4}
\mathscr{L}_\phi w = f \mathscr{L}_\phi ( |\nabla f|^2) +  |\nabla f|^2 \mathscr{L}_\phi f + 2\nabla f\nabla ( |\nabla f|^2).
\end{align}
Putting \eqref{e3.3} and \eqref{e3.4} together gives
\begin{align}\label{e3.5}
(\mathscr{L}_\phi -\partial_t)w  =&   2f|\nabla^2f|^2 + 2f \langle\nabla f, \nabla(\mathscr{L}_\phi-\partial_t)f\rangle + 2f\mathscr{R}ic_\phi(\nabla f,\nabla f) \nonumber \\
& + |\nabla f|^2 (\mathscr{L}_\phi-\partial_t)f + 2\nabla f\nabla(|\nabla f|^2).
\end{align}
Referring to \eqref{e3.1} which simply reads $(\mathscr{L}_\phi -\partial_t)f =  \varepsilon^{-1}(\varepsilon-1)|\nabla f|^2/f - \varepsilon f H(f^{1/\varepsilon})/f^{1/\varepsilon}$ with  $\widehat{H}(f)=H(f^{1/\varepsilon})/f^{1/\varepsilon}$, we can compute
\begin{align*}
\nabla (\mathscr{L}_\phi-\partial_t)f  =  \varepsilon^{-1}(\varepsilon-1) \nabla(|\nabla f|^2)/f - \varepsilon^{-1}(\varepsilon-1)\nabla f|\nabla f|^2/f^2 -\varepsilon \nabla (f\widehat{H}(f))
\end{align*}
so that
\begin{align}\label{e3.6}
2f&  \langle\nabla f, \nabla(\mathscr{L}_\phi-\partial_t)f\rangle\nonumber \\
& = 2   \varepsilon^{-1}(\varepsilon-1) \nabla f \nabla(|\nabla f|^2)- 2  \varepsilon^{-1}(\varepsilon-1)|\nabla f|^4/f - 2\varepsilon f\nabla f \nabla (f\widehat{H}(f))
\end{align}
and
\begin{align}\label{e3.7}
|\nabla f|^2 (\mathscr{L}_\phi-\partial_t)f = \varepsilon^{-1}(\varepsilon-1)|\nabla f|^4/f - \varepsilon|\nabla f|^2f \widehat{H}(f).
\end{align}
Substituting  \eqref{e3.6} and \eqref{e3.7} into \eqref{e3.5} gives the desired result.

\qed

\begin{lemma}\label{lem3.2}
With the assumption of Lemma \ref{lem3.1} we have 
\begin{align}\label{e3.8}
(\mathscr{L}_\phi -\partial_t)w \ge & 2 \varepsilon^{-1}(\varepsilon-1)\frac{1}{f}\langle\nabla f, \nabla w\rangle - \varepsilon^{-1}(5\varepsilon-3)\frac{1}{f^3}w^2  \nonumber\\
& - 2(m-1)\kappa w -  \Gamma^\varepsilon[u] w,
\end{align}
where $\Gamma^\varepsilon[u] = [2H'(u) + (3\varepsilon -2)\widehat{H}(u)]$.
\end{lemma}

\proof
Let us put into consideration the following identity 
$$0\le 2\left|\sqrt{f}\ \nabla^2 f+\frac{\nabla f\otimes \nabla f}{\sqrt{f}} \right|^2 = 2f|\nabla^2 f|^2 + 2\nabla f\nabla |\nabla f|^2 + \frac{2|\nabla f|^4}{f}$$
and the weighted $m$-Ricci tensor bound $\mathscr{R}ic^m_\phi\ge -(m-1)\kappa$ we have from Lemma \ref{lem3.1} as follows
\begin{align}\label{e3.9}
(\mathscr{L}_\phi -\partial_t)w  &= 2\left|\sqrt{f}\ \nabla^2 f+\frac{\nabla f\otimes \nabla f}{\sqrt{f}} \right|^2  + 2\varepsilon^{-1}(\varepsilon-1)\nabla f\nabla(|\nabla f|^2)\nonumber\\
&  \ \ \ \ - (\varepsilon^{-1}(\varepsilon-1)+2)\frac{|\nabla f|^4}{f} + 2f\mathscr{R}ic^m_\phi(\nabla f,\nabla f)  + \frac{2f\langle\nabla \phi,\nabla f\rangle^2}{m-n}\nonumber\\
&  \ \ \ \   -\varepsilon|\nabla f|^2f\widehat{H}(f)- 2\epsilon f \nabla f \nabla (f\widehat{H}(f)) \nonumber\\
& \ge 2\varepsilon^{-1}(\varepsilon-1)\nabla f\nabla(|\nabla f|^2) - \varepsilon^{-1}(3\varepsilon-1)\frac{|\nabla f|^4}{f}  -2(m-1)\kappa f|\nabla f|^2 \nonumber\\
&   \ \ \ \  -\varepsilon|\nabla f|^2f\widehat{H}(f)- 2\epsilon f \nabla f \nabla (f\widehat{H}(f)).
\end{align}
Taking into account the following identity
$$\nabla f\nabla (f |\nabla f|^2) = |\nabla f|^4 + f \nabla f \nabla(|\nabla f|^2)$$
which implies
\begin{align}\label{e3.10}
 \nabla f \nabla(|\nabla f|^2) = ( \nabla f/f)\nabla (f |\nabla f|^2)  - |\nabla f|^4/f
\end{align}
and upon substituting \eqref{e3.9} into \eqref{e3.10} we have 
\begin{align}\label{e3.11}
(\mathscr{L}_\phi -\partial_t)w \ge & 2 \varepsilon^{-1}(\varepsilon-1)\frac{1}{f}\langle\nabla f, \nabla w\rangle - \varepsilon^{-1}(5\varepsilon-3)\frac{1}{f^3}w^2 - 2(m-1)\kappa w  \nonumber\\
&  -\varepsilon|\nabla f|^2f\widehat{H}(f)- 2\epsilon f \nabla f \nabla (f\widehat{H}(f)).
\end{align}
The  explicit value of the last two terms of \eqref{e3.11} can be computed. Recall that  $\widehat{H}(f)=H(f^{1/\varepsilon})/f^{1/\varepsilon}$. We have
\begin{align*}
\nabla (f\widehat{H}(f)) & = \nabla (H(f^{1/\varepsilon})/f^{1/\varepsilon-1})\\
& = \varepsilon^{-1}\nabla f H'(u) + \varepsilon^{-1}(\varepsilon-1)\nabla f \widehat{H}(u)
\end{align*}
so that
\begin{align*}
2\epsilon f \nabla f \nabla (f\widehat{H}(f)) & = 2f|\nabla f|^2 H'(u) + 2(\varepsilon-1)|f\nabla f|^2 \widehat{H}(u)\\
& = 2[H'(u)+ (\varepsilon-1)\widehat{H}(u)]w.
\end{align*}
Therefore we have
\begin{align*}
-\varepsilon|\nabla f|^2f\widehat{H}(f)& - 2\epsilon f \nabla f \nabla (f\widehat{H}(f))  = -2 \left[H'(u)+ \frac{(3\varepsilon-2)}{2}\widehat{H}(u)]\right]w,
\end{align*}
$\widehat{H}(u) = H(u)/u=H(f^{1/\varepsilon})/f^{1/\varepsilon} = \widehat{H}(f)$. Replacing this into \eqref{e3.11} we finally arrive at the required inequality \eqref{e3.8}.

\qed

\subsection*{Proof of Theorem \ref{thm21}}
With recourse to the procedure adopted in the proof of Theorem \ref{thm11} with the same cut-off function $\psi$ described in Lemma \ref{lem23}, we can write using Lemma \ref{lem3.2} that 
\begin{align*}
(\mathscr{L}_\phi -\partial_t)(\psi w) \ge&  2\varepsilon^{-1}(\varepsilon-1)\frac{1}{f}\left[\langle f,\nabla(\psi w)\rangle - w\langle\nabla f,\nabla\psi\rangle \right] \\
& + 2\left[\frac{\nabla \psi}{\psi}\nabla(\psi w) -\frac{|\nabla \psi|^2}{\psi} w \right] - \varepsilon^{-1}(5\varepsilon-3)\frac{1}{f^3}\psi w^2 \\
& + w[\mathscr{L}_\phi-\partial_t -2(m-1)\kappa]\psi -  \Gamma^\varepsilon[u]\psi w.
\end{align*}
Clearly $\psi w$ is positive and attain its maximum at a point in $\mathcal{Q}_{R/2,T/2}(\partial M)$ say $(x_0,t_0)$. Our argument will henceforth split into two cases as in the proof of Theorem \ref{thm11}.

\noindent
{\bf Case 1:} If $x_0 \not\in \partial M$, we may assume withot loss of generality that $x_0\not\in \text{cut}(\partial M)$ by Calabi's argument. Then at the maximum point $(x_0,t_0)$ we have $\nabla(\psi w)= 0$, $\partial_t(\psi w)\ge 0$,  $\mathscr{L}_\phi(\psi w)\le 0$ and then
\begin{align*}
0\ge & - 2\varepsilon^{-1}(\varepsilon-1)\frac{1}{f}\langle\nabla f,\nabla\psi\rangle w - 2\frac{|\nabla \psi|^2}{\psi } w - \varepsilon^{-1}(5\varepsilon-3)\frac{1}{f^3}\psi w^2 \\
& + w[\mathscr{L}_\phi-\partial_t -2(m-1)\kappa]\psi -  \Gamma^\varepsilon[u]\psi w
\end{align*}
which can be rearranged to give 
\begin{align}\label{e3.12}
- \varepsilon^{-1}(5\varepsilon-3) \psi w^2 \le &  2\varepsilon^{-1}(\varepsilon-1)f^2\langle\nabla f,\nabla\psi\rangle w + 2\frac{|\nabla \psi|^2}{\psi }f^3 w \nonumber \\
 & - f^3 w[\mathscr{L}_\phi-\partial_t -2(m-1)\kappa]\psi + f^3 \Gamma^\varepsilon[u]\psi w.
\end{align}
We can now proceed to estimate each term on the right hand side of \eqref{e3.12} as before. This procedure follows directly from the applications of Cauchy-Schwarz and Young's inequalities and properties of $\psi$ as stated in Lemma \ref{lem23}. To this end we estimate the first and second term on the right hand side of \eqref{e3.12} respectively as follows:

\begin{align}\label{e3.13}
2 \frac{(\varepsilon-1)}{\varepsilon}f^2\langle\nabla f,\nabla\psi\rangle w & \le 2 \frac{(1-\varepsilon)}{\varepsilon}\psi^{3/4} w^{3/2} f^{3/2} \frac{|\nabla \psi|}{\psi^{3/4}}\nonumber\\
& \le \frac{(1-\varepsilon)}{\varepsilon}\psi w^2 +C_\varepsilon \frac{|\nabla \psi|^4}{\psi^{3}}(\sup_{\mathcal{Q}_{R,T}} f)^6\nonumber\\
&\le   \frac{(1-\varepsilon)}{6\varepsilon}\psi w^2 + \frac{C_\varepsilon}{R^4}(\sup_{\mathcal{Q}_{R,T}} f)^6
\end{align}
and
\begin{align}\label{e3.14}
2\frac{|\nabla \psi|^2}{\psi}f^3w  & \le 2 \psi^{1/2}w \frac{|\nabla \psi|^2}{\psi^{3/2}}(\sup_{\mathcal{Q}_{R,T}} f)^3\nonumber\\
& \le \frac{(1-\varepsilon)}{6\varepsilon}\psi w^2 +C_\varepsilon \left(\frac{|\nabla \psi|^2}{\psi^{3/2}}\right)^2(\sup_{\mathcal{Q}_{R,T}} f)^6\nonumber\\
& \le \frac{(1-\varepsilon)}{6\varepsilon}\psi w^2  +  \frac{C_\varepsilon}{R^4}(\sup_{\mathcal{Q}_{R,T}} f)^6
\end{align}
at $(x_0,t_0)$. 

\noindent
Similarly to \eqref{A9} - \eqref{A11} with little modification we have  
for the term with time derivative of $\psi$:
\begin{align}\label{e3.15}
f^3w\psi_t & \le \psi^{\frac{1}{2}} w \frac{|\psi_t|}{\psi^{\frac{1}{2}}}\Big(\sup_{\mathcal{Q}_{R,T}} f\Big)^3 \nonumber\\
& \le   \frac{(1-\varepsilon)}{6\varepsilon}\psi w^2  +C_\varepsilon\frac{|\psi_t|^2}{\psi} \Big(\sup_{\mathcal{Q}_{R,T}} f\Big)^6\nonumber\\
& \le \frac{(1-\varepsilon)}{6\varepsilon}\psi w^2  + \frac{C_\varepsilon}{T^2} \Big(\sup_{\mathcal{Q}_{R,T}} f\Big)^6.
\end{align}
For the term with $(m-1)\kappa$ we have
\begin{align}\label{e3.16}
2f^3 (m-1)\kappa\psi w & \le  2(m-1)\kappa  (\psi^{\frac{1}{2}}w) \psi^{\frac{1}{2}}\Big(\sup_{\mathcal{Q}_{R,T}} f\Big)^3  \nonumber\\
& \le  \frac{(1-\varepsilon)}{6\varepsilon}\psi w^2 + C_\varepsilon m^2\kappa^2 \Big(\sup_{\mathcal{Q}_{R,T}} f\Big)^6.
\end{align}
For the term involving the weighted Laplacian of $\psi$ we have
\begin{align*}
-f^3 w \mathscr{L}_\phi\psi  & =-  f^3  \left[ \psi_\varrho \mathscr{L}_\phi\varrho + \psi_{\varrho \varrho} |\nabla \varrho|^2 \right]w \\
&\le \left[| \psi_\varrho| ((n-1)\kappa R+\ell) + \psi_{\varrho \varrho} |\nabla \varrho|^2 \right]w \Big(\sup_{\mathcal{Q}_{R,T}} f\Big)^3
\end{align*}
which through similar step as in \eqref{A11} yields the bound
\begin{align}\label{e3.17}
-f^3 w \mathscr{L}_\phi\psi & \le \frac{(1-\varepsilon)}{6\varepsilon}\psi w^2  +  \left( \frac{C_\varepsilon}{R^4}  + C_\varepsilon \kappa^2 + C_\varepsilon \ell^4\right)\Big(\sup_{\mathcal{Q}_{R,T}} f \Big)^6.
\end{align}
 since $\mathscr{R}ic_\phi^m \ge -(m-1)\kappa$ and $\mathscr{H}_\phi \ge -\ell$. 
Combining \eqref{e3.15} - \eqref{e3.17} we have the bound on the third term on the right hand side of \eqref{e3.12} as follows
\begin{align}\label{e3.18}
- f^3 w[\mathscr{L}_\phi & -\partial_t -2(m-1)\kappa]\psi \nonumber \\
 & \le \frac{(1-\varepsilon)}{2\varepsilon}\psi w^2  +  C_\varepsilon \left( \frac{1}{R^4}  + \kappa^2 +  \ell^4  +\frac{1}{T^2}\right)\Big(\sup_{\mathcal{Q}_{R,T}} f \Big)^6.
\end{align}
at $(x_0,t_0)$.

\noindent
For the last term we compute at $(x_0,t_0)$
\begin{align}
 f^3 \Gamma^\varepsilon[u]\psi w & \le \psi^{\frac{1}{2}}w  [ \Gamma^\varepsilon[u]]_+ \psi^{\frac{1}{2}}\Big(\sup_{\mathcal{Q}_{R,T}} f\Big)^3  \nonumber\\
& \le  \frac{(1-\varepsilon)}{6\varepsilon}\psi w^2 + C_\varepsilon [ \Gamma^\varepsilon[u]]_+^2 \Big(\sup_{\mathcal{Q}_{R,T}} f\Big)^6,
\end{align}
where $[ \Gamma^\varepsilon[u]]_+ = [2H'(u)+(3\varepsilon-2)\widehat{H}(u)]_+$.

\noindent
Putting these estimates back into \eqref{e3.12} and rearranging we obtain
\begin{align*}
\varepsilon^{-1}(3-5\varepsilon) \psi w^2 \le &  2\varepsilon^{-1}(1-\varepsilon)\psi w^2 +  C_\varepsilon \left( \frac{1}{R^4}  +\frac{1}{T^2} + \kappa^2 +  \ell^4   +  [ \Gamma^\varepsilon[u]]_+^2 \right)\Big(\sup_{\mathcal{Q}_{R,T}} f \Big)^6
\end{align*}
at $(x_0,t_0)$.  It therefore follows that 
\begin{align*}
2\varepsilon^{-1}(1-2\varepsilon) \psi w^2 \le &    C_\varepsilon \left( \frac{1}{R^4}  +\frac{1}{T^2} + \kappa^2 +  \ell^4   + [ \Gamma^\varepsilon[u]]_+^2 \right)\Big(\sup_{\mathcal{Q}_{R,T}} f \Big)^6
\end{align*}
by restricting $\varepsilon$ to the range $0<\varepsilon <1/2$ at $(x_0,t_0)$.  Hence
\begin{align*}
\psi w^2 \le C_\varepsilon  \left( \frac{1}{R^4}  +\frac{1}{T^2} + \kappa^2 +  \ell^4   +  [ \Gamma^\varepsilon[u]]_+^2 \right)\Big(\sup_{\mathcal{Q}_{R,T}} f \Big)^6
\end{align*}
at $(x_0,t_0)$
Reverting to the maximal characterization of the point $(x_0,t_0)$ we have for all $(x,t)\in \mathcal{Q}_{R/2,T/2}(\partial M)$
\begin{align*}
\psi^2w^2(x,t) & \le \psi^2w^2(x_0,t_0)\le \psi w^2(x_0,t_0).
\end{align*}
Since $\psi(x,t)\equiv 1$ on $\mathcal{Q}_{R/2,T}(\partial M)$ and by definition  of $w$, i.e., $ w = f |\nabla f|^2$  with $f=u^\varepsilon$, $\varepsilon \in (0, 1/2)$ we obtain
\begin{align*}
 (f|\nabla f|^2)^2 \le C_\varepsilon \left( \frac{1}{R^4}  +\frac{1}{T^2} + \kappa^2 +  \ell^4   +  [ \Gamma^\varepsilon[u]]_+^2 \right)\Big(\sup_{\mathcal{Q}_{R,T}} f \Big)^6.
\end{align*}
Noting that 
$f|\nabla f|^2 = \varepsilon^2 f^3|\nabla u|^2/u^2$ and $(\sup\{f\}) = (\sup \{u^\varepsilon\})$. Then we arrive at 
\begin{align*}
\frac{|\nabla u|^2} {u^{2-3\varepsilon}} \le C_\varepsilon \left( \frac{1}{R^2}  +\frac{1}{T} +L  +  \sup_{\mathcal{Q}_{R,T}}  [ \Gamma^\varepsilon[u]]_+ \right)\Big(\sup_{\mathcal{Q}_{R,T}} u^\varepsilon\Big)^3
\end{align*}
for all $(x,t)\in \mathcal{Q}_{R/2,T/2}(\partial M)$, where $L= (\kappa^2 +  \ell^4)^{1/2}$. In conclusion, rearranging gives the desired estimates.
 
\noindent
{\bf Case 2:}
If $x_0 \in \partial M$, we will show that the gradient estimates \eqref{ethm21} holds on $\partial M \times [0, T]$.  Indeed, at the maximum point $(x_1,t_1)  \in \partial M \times [0, T]$ of $\psi w$ we have $(\psi w)_\nu \ge 0$ which implies $\psi_\nu w +\psi w_\nu = \psi w_\nu \ge 0$, and in particular 
$$w_\nu \ge 0.$$
Since $w=f|\nabla f|^2$, and $f=u^\varepsilon$ with $\varepsilon \in (0,1)$ satisfies the Dirichlet boundary condition we have 
\begin{align}\label{e3.21}
0\le w_\nu = (f|\nabla f|^2)_\nu = f_\nu |\nabla f|^2 + f(|\nabla f|^2)_\nu.
\end{align}
Noting that $|\nabla f|=f_\nu = \varepsilon u^\varepsilon u_\nu/u\ge 0$ since $u>0$ and $u_\nu\ge 0$.  Now referring to the Reilly-type formula formula \eqref{Rei} then \eqref{e3.21} gives
\begin{align}\label{e3.22}
0\le f_\nu f_\nu^2 + 2ff_\nu \left[\mathscr{L}_\phi f - f_\nu \mathscr{H}_\phi\right].
\end{align}
To keep the analysis simple, we write (c.f. \eqref{e3.1} and note that $\widehat{H}(f^{1/\varepsilon}) = H(f^{1/\varepsilon})/f^{1/\varepsilon}$)
\begin{align*}
\mathscr{L}_\phi f & = \varepsilon^{-1}(\varepsilon-1)|\nabla f|^2/f + f_t - \varepsilon f\widehat{H}(f^{1/\varepsilon})\\ 
& = \varepsilon^{-1}(\varepsilon-1)f_\nu^2/f + f_t - \varepsilon f\widehat{H}(f^{1/\varepsilon}).
\end{align*}
Noting that 
$f^2_\nu =  w/f$ and $f_\nu f_\nu^2 = (w/f)^{3/2}.$ It the follows by plugging in these data into \eqref{e3.22} that 
\begin{align}\label{e3.23}
0\le \frac{w^{3/2}}{f^{3/2}} + \frac{2 w^{1/2}}{f^{-1/2}}\left[\frac{\varepsilon^{-1}(\varepsilon-1)w}{f^2} +f_t -\varepsilon f\widehat{H}(f^{1/\varepsilon}) -  \frac{w^{1/2}}{f^{1/2}} \mathscr{H}_\phi  \right]
\end{align}
at $(x_0,t_0)$. Within the consideration of the condition that $u_t\le H(u)$ which implies $f_t -\varepsilon f\widehat{H}(f^{1/\varepsilon})\le 0$ over $\partial M\times [0,T]$ we see that \eqref{e3.23} is reduced to
\begin{align}\label{e3.24}
-[1+2\varepsilon^{-1}(\varepsilon-1)]\frac{w^{3/2}}{f^{3/2}} + 2w  \mathscr{H}_\phi \le 0,
\end{align}
Since $\varepsilon \in (0,1/2)$ then $-[1+2\varepsilon^{-1}(\varepsilon-1)] = 2\varepsilon^{-1}(1-\varepsilon)-1>0$. Clearly \eqref{e3.24} can be solved to give
$$w(x_0,t_0) = 0 \ \ \text{or}\ \ \ w^{1/2}(x_0,t_0) \le \frac{2}{2\varepsilon^{-1}(1-\varepsilon)-1} \ell (\sup f)^{3/2}$$
since $ \mathscr{H}_\phi  \ge -\ell$ by the theorem hypothesis. It therefore means that 
$$\psi w=0 \ \ \text{or} \ \ \psi w(x_0,t_0) \le \frac{4}{(2\varepsilon^{-1}(1-\varepsilon)-1)^2}\ell^2 (\sup f)^3$$
on $\mathcal{Q}_{R,T}(\partial M)$. The rest of the proof is exactly the same as that of Theorem \ref{thm11}. Thus 
\begin{align*}
f|\nabla f|^2(x,t) =w(x,t)=\psi w (x,t)\le \psi w(x_0,t_0) \le  \frac{4}{(2\varepsilon^{-1}(1-\varepsilon)-1)^2}\ell^2 (\sup f)^3,
\end{align*}
which also implies the conclusion by using the relation $f = u^\varepsilon$  
with $\varepsilon$ restricted to the range $(0,1/2)$.

\section{Li-Yau type gradient estimates}\label{sec5}

\begin{theorem}\label{thm13}
Let $M=(M,g,e^{-\phi}dv)$ be an $n$-dimensional complete smooth metric measure space with compact boundary satisfying the bounds
\begin{align*}
\mathscr{R}ic_\phi^m \ge -(m-1)k \ \ \text{and}\ \ \mathscr{H}_\phi \ge - l
\end{align*}
for some constants $k,l \ge 0$.  Let $u=u(x,t)$ be a smooth positive solution to \eqref{11} in $\mathcal{Q}_{2R,T}(\partial M)$. Suppose further that $u$ satisfies the Dirichlet boundary condition (i.e., $u(\cdot,t)|_{\partial M}$ is constant for each time $t \in (-\infty,+\infty)$),
\begin{align*}
u_\nu\ge 0\ \ \text{and}\ \ u_t \le H(u)
\end{align*}
over $\partial M \times [T_0-T, T_0]$.   Then have for all $(x,t)\in \mathcal{Q}_{R,T}(\partial M)$, $t>0$
\begin{align}\label{eqthm3}
\frac{|\nabla u|^2}{u^2} &+\delta\frac{H(u)}{u} -\delta \frac{u_t}{u}
 \le \frac{m\delta^2}{2(1-\epsilon)}\Big(\frac{1}{t}+ \frac{k}{\delta-1}+ \widetilde{\Phi} +\mu^H\Big) +  \frac{\delta}{p(\delta-1)} \sqrt{ \frac{m}{2(1-\epsilon)}}\gamma^H
\end{align}
on $\partial M \times [T_0-T, T_0]$.  Here  $\epsilon \in (0,1)$, $\delta>1$ and  $p>0$ such that 
\begin{align}\label{hyp5}
\frac{2(1-\epsilon)(\delta-1)}{m\delta p}\ge \frac{1}{\epsilon}-1+\frac{(1-\epsilon)(1-\delta)^2}{4},
\end{align}
and
$$\widetilde{\Phi} = \frac{1}{R^2}\Big[\Big(\frac{m\delta^2}{4(1-\epsilon)(\delta p+\delta-1)}\Big)C_1^2 +C_2 +  C_1 [(m-1)kR+l]R\Big]\Big]$$
for some constants positive $C_1$ and $C_2$.
Furthermore,
$$\mu^H : =   \sup_{\mathcal{Q}_{2R,T} } \{[H_u]_+\}  \ \ \text{and}\ \ \gamma^H := \sup_{\mathcal{Q}_{2R,T} }  \{[\Omega^H]_+\},$$ 
where $[H_u]_+ = [H'(u)-H(u)/u]_+$ and $[\Omega^H]_+ = [-[H'(u)-H(u)/u-(\alpha/p^2)uH''(u)]]_+$.
\end{theorem}

\subsection{Evolution equations}
Here in this section we consider the nonlinear heat-type equation
\begin{align}\label{e4.1}
u_t = \mathscr{L}_\phi u + H(u).
\end{align}
Let
$$w=u^{-p},$$
where $p>0$ is to be determined later. By direct computation as before we obtain
\begin{align*}
w_t  = -pu^{-p-1}u_t, \ \ \nabla w  = -pu^{-p-1}\nabla u,\ \ |\nabla w|^2/w^2 =p^2|\nabla u|^2/u^2
\end{align*}
and
\begin{align*}
\Delta w = -pu^{-p-1}\Delta u + p(p-1)u^{-p-2}|\nabla u|^2.
\end{align*}
Clearly by the definition of $\mathscr{L}_\phi w= \Delta w -\nabla\phi \nabla w$ we have
\begin{align*}
\mathscr{L}_\phi w & = -pu^{-p-1}\mathscr{L}_\phi u + p(p-1)u^{-p-2}|\nabla u|^2\\
& =  -pu^{-p-1}(u_t-H(u)) + ((p+1)/p)|\nabla w|^2/w.
\end{align*}
Then we have proved the following lemma.

\begin{lemma}\label{lem41}
Suppose $u>0$ is a smooth solution of \eqref{e4.1}. The the function $w=u^{-p}$, $p>0$,  satisfies the following equation
\begin{align}\label{e4.2}
(\mathscr{L}_\phi -\partial_t)w & = \frac{P+1}{p}\frac{|\nabla w|^2}{w} + pw\widehat{H}(w),
\end{align}
where $\widehat{H}(w) = H(w^{-1/p})/w^{-1/p} = H(u)/u$.
\end{lemma}
For the purpose of Li-Yau type gradient estimate to be established in this section, we introduce the following Harnack quantity
\begin{align}\label{e4.3}
F(x,t) =F_0(x,t) +\beta F_1(x,t).
\end{align}
For computational purpose, we also introduce $F_0$ and $F_1$ in the following quantities
\begin{align*}
F_0(x,t)=  \frac{|\nabla w|^2}{w^2} + \alpha\widehat{H}(w)  \ \ \text{and} \ \ \ F_1(x,t) = \frac{w_t}{w},
\end{align*}
where $\alpha,\beta>0$ are constants to be chosen later.  Before proceeding to prove the desired Li-Yau type gradient estimate we first prove some computational lemmas.

\begin{lemma}\label{lem42}
Suppose $u>0$ is a smooth solution of \eqref{e4.1} and $w$ satisfies \eqref{e4.2}. Then 
\begin{align}\label{e4.4}
(\mathscr{L}_\phi -\partial_t)F_1(x,t) &= \frac{2}{p}\nabla F_1(x,t) \nabla \log w(x,t) + p\widehat{H}_t.
\end{align}
\begin{align}\label{e4.5}
(\mathscr{L}_\phi -\partial_t)F_0(x,t) &= \frac{\mathscr{L}_\phi|\nabla w|^2}{w^2} +2\Big(3-\frac{p+1}{p}\Big)\frac{|\nabla w|^4}{w^4} - \frac{8\langle \nabla^2w,\nabla w\otimes \nabla w \rangle}{w^3}\nonumber\\
& \hspace{1cm} - \frac{2\nabla w\nabla w_t}{w^2} + \alpha \Big(\mathscr{L}_\phi -\partial_t - \frac{2p}{\alpha}\frac{|\nabla w|^2}{w^2}\Big)\widehat{H}.
\end{align}
\end{lemma}
\proof
A simple computation gives 
\begin{align*}
\partial_tF_1(x,t) & =\frac{1}{w^2}(ww_{tt}-w^2_t)\\
\nabla F_1(x,t) & = \frac{\nabla w_t}{w} - \frac{w_t\nabla w}{w^2}\\
\mathscr{L}_\phi F_1(x,t) & = e^{\phi}\text{div}(e^{-\phi}F_1) = \text{div}\nabla F_1 - \nabla\phi\nabla F_1\\
& = \frac{1}{w^2} (w \mathscr{L}_\phi  w_t - 2\nabla w\nabla w_t -w_t\mathscr{L}_\phi  w + 2w_t |\nabla w|^2/w).
\end{align*}
Combining the expressions for $\partial_t F_1$ and $\mathscr{L}_\phi F_1$ and rearranging terms of the resulting expression, we have
\begin{align}\label{e4.6}
(\mathscr{L}_\phi -\partial_t)F_1(x,t) =  \frac{1}{w}(\mathscr{L}_\phi w-w_t)_t -\frac{w_t}{w^2}(\mathscr{L}_\phi w-w_t) -\frac{2\nabla w\nabla w_t}{w^2}+ 2w_t\frac{|\nabla w|^2}{w^3}.
\end{align}
Substituting \eqref{e4.2} into \eqref{e4.6} and simplifying gives 
\begin{align*}
(\mathscr{L}_\phi -\partial_t)F_1(x,t) & = \frac{1}{w}\Big(\frac{P+1}{p}\frac{|\nabla w|^2}{w} + pw\widehat{H}\Big)_t - \frac{w_t}{w^2}\Big(\frac{P+1}{p}\frac{|\nabla w|^2}{w} + pw\widehat{H}\Big)\\
& \hspace{1cm}  -\frac{2\nabla w\nabla w_t}{w^2}+ 2w_t\frac{|\nabla w|^2}{w^3}\nonumber\\
& = \frac{p+1}{p}\Big(\frac{2\nabla w\nabla w_t}{w^2} - \frac{2 |\nabla w|^2  w_t}{w^3} \Big)  + p\widehat{H}_t -  \frac{p+1}{p} \frac{|\nabla w|^2  w_t}{w^3}\nonumber\\
& \hspace{1cm}  -\frac{2\nabla w\nabla w_t}{w^2}+ 2w_t\frac{|\nabla w|^2}{w^3}\nonumber\\
& = \frac{2}{p}\nabla F_1\nabla \log w + p\widehat{H}_t.
\end{align*}
Similarly
\begin{align*}
\partial_t F_0(x,t) = \frac{2\nabla w\nabla w_t}{w^2} - \frac{2|\nabla w|^2w_t}{w^3} +  \alpha \widehat{H}_t.
\end{align*}

\begin{align*}
\mathscr{L}_\phi  F_0(x,t) & = \frac{\mathscr{L}_\phi |\nabla w|^2}{w^2} +|\nabla w|^2 \mathscr{L}_\phi \left(\frac{1}{w^2}\right) + 2\Big\langle\nabla |\nabla w|^2, \nabla \left(\frac{1}{w^2}\right) \Big\rangle + \alpha \widehat{H}_t\\
& = \frac{\mathscr{L}_\phi |\nabla w|^2}{w^2}  + \frac{6|\nabla w|^4}{w^4} - \frac{2|\nabla w|^2  \mathscr{L}_\phi w}{w^3}-  \frac{8\langle \nabla^2w,\nabla w\otimes \nabla w \rangle}{w^3} + \alpha \widehat{H}_t
\end{align*}
and then 
\begin{align*}
(\mathscr{L}_\phi -\partial_t)F_0(x,t) &= \frac{\mathscr{L}_\phi|\nabla w|^2}{w^2} + \frac{6|\nabla w|^4}{w^4} - \frac{8\langle \nabla^2w,\nabla w\otimes \nabla w \rangle}{w^3}  - \frac{2 \nabla w\nabla w_t}{w^2}  \nonumber\\
& \hspace{1cm} - \frac{2|\nabla w|^2}{w^3}(\mathscr{L}_\phi w-w_t)  + \alpha (\mathscr{L}_\phi -\partial_t)\widehat{H}.
\end{align*}
Substituting \eqref{e4.2} into the last expression gives \eqref{e4.5}. This complete the proof.

\qed

\begin{lemma}\label{lem43}
With the assumption of Lemma \ref{lem42} we have for $\epsilon \in (0,1)$ 
\begin{align}\label{e4.7}
(\mathscr{L}_\phi -\partial_t)F_0(x,t) & \ge \frac{2(1-\epsilon)}{m}\Big(\frac{\mathscr{L}_\phi w}{w}\Big)^2 - 2\Big(\frac{1}{\epsilon}-1\Big)\frac{|\nabla w|^2}{w^4} + \frac{2}{w^2}\mathscr{R}ic_\phi^m(\nabla w, \nabla w) \nonumber\\
& \hspace{1cm} + \frac{2}{p}\nabla F_0\nabla \log w + 2\Big(p-\frac{\alpha}{p}\Big) \frac{\nabla w}{w}\nabla \widehat{H} + \alpha (\mathscr{L}_\phi -\partial_t)\widehat{H}.
\end{align}
\end{lemma}

\proof
First, by applying the weighted Bochner formula \eqref{e1} we write the evolution equation \eqref{e4.5} for $F_0(x,t)$ to have 
\begin{align}\label{e4.8}
(\mathscr{L}_\phi -\partial_t)F_0(x,t) & \ge \Big(\frac{2|\nabla^2w|^2}{w^2} + \frac{6|\nabla w|^4}{w^4} - \frac{8\langle \nabla^2w,\nabla w\otimes \nabla w \rangle}{w^3}\Big)  + \frac{2\mathscr{R}ic_\phi(\nabla w, \nabla w)}{w^2}\nonumber\\
& \hspace{1cm} + \frac{2\langle \nabla w, \nabla(\mathscr{L}_\phi w-w_t)\rangle}{w^2} - \frac{2(p+1)}{p}\frac{|\nabla w|^4}{w^4}\nonumber \\
& \hspace{1cm} + \alpha \Big(\mathscr{L}_\phi -\partial_t - \frac{2p}{\alpha}\frac{|\nabla w|^2}{w^2}\Big)\widehat{H}.
\end{align}
We will combine the evolution equations for $F_0(x,t)$ and $F_1(x,t)$ as presented in Lemma \ref{lem42}, but before we proceed lets do further analysis on \eqref{e4.8}. To do this we find some lower bounds as follows. By Young's inequality for $\epsilon \in (0,1)$ we have
\begin{align*}
\frac{4\langle\nabla^2w, \nabla w\otimes\nabla w\rangle}{w^3}\le \frac{4|\nabla^2 w|}{w}\frac{|\nabla w\otimes\nabla w|}{w^2} \le \frac{2\epsilon|\nabla^2 w|^2}{w^2}+\frac{2}{\epsilon}\frac{|\nabla w|^4}{w^2}
\end{align*}
and then
\begin{align}\label{e4.9}
& \frac{2|\nabla^2w|^2}{w^2}  + \frac{6|\nabla w|^4}{w^4}  - \frac{8\langle \nabla^2w,\nabla w\otimes \nabla w \rangle}{w^3}\\
& \ge 2(1-\epsilon)\frac{|\nabla^2w|^2}{w^2} -4\Big(\frac{\langle \nabla^2w,\nabla w\otimes \nabla w \rangle}{w^3} -   \frac{|\nabla w|^4}{w^4} \Big) - 2\Big(\frac{1}{\epsilon}-1\Big) \frac{|\nabla w|^4}{w^4}.\nonumber
\end{align}
Noting also by Newton inequality \eqref{NI} (also check for \eqref{ee1}) we can write
\begin{align*}
2(1-\epsilon)\frac{|\nabla^2w|^2}{w^2} &  \ge \frac{2(1-\epsilon)}{w^2}\Big(\frac{1}{m}(\mathscr{L}_\phi w)^2  - \frac{\langle\nabla\phi,\nabla w\rangle^2}{(m-n)}\Big)\\
& \ge  \frac{2(1-\epsilon)}{m} \Big(\frac{\mathscr{L}_\phi w}{w}\Big)^2 -  \frac{2\langle\nabla\phi,\nabla w\rangle^2}{w^2(m-n)}
\end{align*}
since $2\epsilon\langle\nabla\phi,\nabla w\rangle^2/(w^2(m-n))\ge 0$   for $0<\epsilon<1$. Thus we have 
\begin{align}\label{e4.10}
2(1-\epsilon)\frac{|\nabla^2 w|^2}{w^2}  + \frac{2\mathscr{R}ic_\phi(\nabla w, \nabla w)}{w^2} \ge \frac{2(1-\epsilon)}{m}\Big(\frac{\mathscr{L}_\phi w}{w}\Big)^2 + \frac{2\mathscr{R}ic_\phi^m(\nabla w, \nabla w)}{w^2}.
\end{align}
In consideration of Lemma \ref{lem41} we compute the term
\begin{align}\label{e4.11}
 \frac{2\langle \nabla w, \nabla(\mathscr{L}_\phi w-w_t)\rangle}{w^2} &= \frac{4(p+1)}{p} \frac{\langle \nabla^2w,\nabla w\otimes \nabla w \rangle}{w^3} -  \frac{2(p+1)}{p}\frac{|\nabla w|^4}{w^4} \nonumber\\
 & \hspace{1cm} + 2p\frac{|\nabla w|^2}{w^2}\widehat{H} + 2p
\frac{\langle\nabla w, \nabla \widehat{H}\rangle}{w}.
\end{align}
Substituting these equations \eqref{e4.9}-\eqref{e4.11} into \eqref{e4.8} we have
\begin{align}\label{e4.12}
(\mathscr{L}_\phi -\partial_t)F_0(x,t) & \ge   \frac{2(1-\epsilon)}{m}\Big(\frac{\mathscr{L}_\phi w}{w}\Big)^2  - 2\Big(\frac{1}{\epsilon}-1\Big) \frac{|\nabla w|^4}{w^4} + \frac{2\mathscr{R}ic_\phi^m(\nabla w, \nabla w)}{w^2}\nonumber\\
& +\frac{4}{p}\Big(\frac{\langle \nabla^2w,\nabla w\otimes \nabla w \rangle}{w^3} -   \frac{|\nabla w|^4}{w^4} \Big) +  2p\frac{\nabla w\nabla \widehat{H}}{w}  \nonumber\\
& +  \alpha \Big(\mathscr{L}_\phi -\partial_t \Big)\widehat{H}.
\end{align}
The middle line on the right hand side of  \eqref{e4.12} (that is,  the fourth and fifth terms on the right hand side of  \eqref{e4.12}) are simplified to get
\begin{align}\label{e4.13}
\frac{4}{p}\Big(\frac{\langle \nabla^2w,\nabla w\otimes \nabla w \rangle}{w^3}  &- \frac{|\nabla w|^4}{w^4}\Big)  +  2p\frac{\nabla w,}{w} \nabla \widehat{H} \nonumber\\
& = \frac{2}{p}\nabla F_0\nabla\log w + 2\Big(p-\frac{\alpha}{p}\Big)\frac{\nabla w\nabla \widehat{H}}{w}. 
\end{align}
Putting \eqref{e4.13} into \eqref{e4.12} gives \eqref{e4.7} which is the desired evolution inequality for $F_0(x,t)$.

\qed

Now that we have evolution formulas for $F_0$ and $F_1$ on our hands the next is to obtain evolution formula for $F$. Recall that $F=F_0+\beta F_1$ for $\beta>0$. Then combining \eqref{e4.4} of Lemma \ref{lem42} and \eqref{e4.7} of Lemma \ref{lem43} choosing $\beta=\alpha/p$ we have
\begin{align}\label{e4.14}
(\mathscr{L}_\phi -\partial_t)F & \ge \frac{2(1-\epsilon)}{m}\Big(\frac{\mathscr{L}_\phi w}{w}\Big)^2 - 2\Big(\frac{1}{\epsilon}-1\Big)\frac{|\nabla w|^4}{w^4} + \frac{2\mathscr{R}ic_\phi^m(\nabla w, \nabla w)}{w^2} \nonumber\\
& \hspace{1cm} + \frac{2}{p}\nabla F\nabla \log w + 2\Big(p-\frac{\alpha}{p}\Big) \frac{\nabla w}{w}\nabla \widehat{H} + \alpha \mathscr{L}_\phi \widehat{H}.
\end{align}
 Clearly by Lemma \ref{lem41} with $\beta=\alpha/p$  we have
 \begin{align*}
(\mathscr{L}_\phi w)/w = [(p+1)/p-p/\alpha](|\nabla w|^2/w^2 )+ (p/\alpha)F.
 \end{align*}
Therefore, letting $\alpha=\delta p^2$ for $\delta>1$ yields
\begin{align}\label{e4.15}
(\mathscr{L}_\phi w)/w = (\delta p+\delta-1)/ (\delta p)(|\nabla w|^2/w^2 ) +  (1/\delta p) F.
\end{align}
Substituting \eqref{e4.15} into \eqref{e4.14} we arrive at

\begin{align}\label{e4.16}
(\mathscr{L}_\phi -\partial_t)F  \ge & \frac{2(1-\epsilon)}{m\delta^2p}F^2 +\Big[  \frac{2(1-\epsilon)}{m}\frac{(\delta p+\delta-1)^2}{\delta^2p^2}- 2\Big(\frac{1}{\epsilon}-1\Big)\Big] \frac{|\nabla w|^4}{w^4}   \nonumber\\
&+ \frac{4(1-\epsilon)}{m}\frac{(\delta p+\delta-1)}{\delta^2p^2}F \frac{|\nabla w|^2}{w^2}
+ \frac{2}{p}\nabla F\nabla \log w \nonumber\\ 
& + \frac{2\mathscr{R}ic_\phi^m(\nabla w, \nabla w)}{w^2} + 2p(1-\delta) \frac{\nabla w}{w}\nabla \widehat{H} + \delta p^2  \mathscr{L}_\phi \widehat{H}.
\end{align}
We can also compute the explicit value of the sum of the last two terms on the right hand side of \eqref{e4.16} using $\widehat{H} =H(u)/u$, $w=u^{-p}$. First note that 
$$\nabla \widehat{H} = \nabla (H(u)/u) = H_u(\nabla u/u) = -(1/p)H_u\nabla w/w,$$
where $H_u:=H'(u)-H(u)/u$.  Using this idea we thus compute the second to the last term of \eqref{e4.16} as 
\begin{align}\label{e4.17}
2p(1-\delta) \frac{\nabla w}{w}\nabla \widehat{H} = -2(1-\delta)H_u |\nabla w|^2/w^2.
\end{align}
Likewise we compute the last term of \eqref{e4.16} as follows
\begin{align*}
\mathscr{L}_\phi \widehat{H} & = e^\phi\text{div}(e^{-\phi}\nabla \widehat{H}) = \nabla\cdot(\nabla \widehat{H})-\nabla\phi\nabla\widehat{H}\\
& = \nabla\cdot(H_u(\nabla u/u))- H_u \nabla\phi\nabla u/u\\
& = (uH''(u)-H_u)|\nabla u|^2/u^2 + \nabla\cdot(\nabla u/u)- H_u \nabla\phi\nabla u/u\\
& = (uH''(u)-2H_u)|\nabla u|^2/u^2 + H_u (\mathscr{L}_\phi u)/u.
\end{align*}
Using $|\nabla u|^2/u^2 = (1/p^2)|\nabla w|^2/w^2$ and $\mathscr{L}_\phi u=u_t-H(u)$ we have 
\begin{align}\label{e4.18}
\alpha \mathscr{L}_\phi \widehat{H}  = (\alpha/p^2)(uH''(u)-2H_u)|\nabla w|^2/w^2 -(\alpha/p)H_u w_t/w -\alpha H_u H(u)/u.
\end{align}
Adding the expressions \eqref{e4.17} and \eqref{e4.18} and using $\alpha=\delta p^2$ and $\beta=\alpha/p$, rearranging terms and observing some cancellation we obtain 
\begin{align*}
\alpha \mathscr{L}_\phi \widehat{H} & + 2p(1-\delta) \frac{\nabla w}{w}\nabla \widehat{H}  = \Omega_H|\nabla w|^2/w^2  - H_u F,
\end{align*}
where $\Omega_H:= -[H_u - (\alpha/p^2) uH''(u)]$ and $H_u:=H'(u)-H(u)/u$.
Combining the above discussion with Bakry-\'Emery Ricci tensor lower bound, $\mathscr{R}ic_\phi^m\ge -(m-k)g$,  $k\ge 0$ we arrive at the final lemma in this section.
\subsection{Fundamental lemma}
\begin{lemma}\label{lem44}
With the assumption of Lemma \ref{lem43}  and if in addition, $\mathscr{R}ic_\phi^m\ge -(m-k)g$,  $k\ge 0$  holds we have 
\begin{align}\label{e4.19}
(\mathscr{L}_\phi -\partial_t)F  \ge & \frac{2(1-\epsilon)}{m\delta^2p}F^2 +\Big[  \frac{2(1-\epsilon)}{m}\frac{(\delta p+\delta-1)^2}{\delta^2p^2}- 2\Big(\frac{1}{\epsilon}-1\Big)\Big] \frac{|\nabla w|^4}{w^4}   \nonumber\\
&+ \frac{4(1-\epsilon)}{m}\frac{(\delta p+\delta-1)}{\delta^2p^2}F \frac{|\nabla w|^2}{w^2}
+ \frac{2}{p}\nabla F\nabla \log w \nonumber\\ 
& - 2(m-1)k\frac{|\nabla w|^2}{w^2} + \Omega^H|\nabla w|^2/w^2  - H_u F,
\end{align}
where $\Omega^H:= -[H_u - (\alpha/p^2) uH''(u)]$ and $H_u:=H'(u)-H(u)/u$.
\end{lemma}

The last lemma will be combined with suitable cut-off function to finally prove Theorem \ref{thm13} in the remaining part of this section. Before we proceed we briefly summarize the properties of such cut-off function as propagated by Li-Yau \cite{[LY86]}.

\subsection{A suitable cut-off function and its consequence}
As in the previous sections, the distance function from the boundary $\varrho_{\partial M} : M\to \mathbb{R}^+$ is defined $\varrho_{\partial M}=d(\cdot,\partial M)$ which is smooth outside the cut locus of the boundary $\text{Cut}(\partial M)$. For a reference point in $M$ we fix $T,R>0$ , and define a cut-off function $\chi(r)$ in the interval $[0,+\infty)$. 
 We  set 
\begin{align}\label{eC}
\psi(\varrho_{\partial M}(x))=\chi\left(\frac{\varrho_{\partial M}(x)}{R}\right)
\end{align}
for $x\in \partial M$ and $t>0$.  The following existence result due to  Li-Yau \cite{[LY86]} states its useful properties to the next localization procedure.

\begin{lemma}\label{lem45}
There exists a smooth non-increasing cut-off-function $\chi(r)$  atleast $C^2([0,+\infty))$ such that 
\begin{enumerate}
\item  $0\leq \chi(r) \leq 1$ for $0\le r\le 1$, $\chi(r)\equiv 1$ for $r\le 1$ and $\chi(r)\equiv 0$ for $r\ge 2$.
\item  $-c_1 \chi^{1/2} \le \chi'(r) \le 0$   and $\chi''(r) \ge -C_2$ for $r\in [0,+\infty)$ and some constants $C_1,C_2>0$.
\item $\chi(r)$ is non-increasing as a radial function of r on $[0,+\infty)$.
\end{enumerate}
\end{lemma}

Since the weighted Ricci curvature tensor is bounded below by $\mathscr{R}ic_\phi^m\ge -(m-1)k$, 
$k\ge 0$ and the weighted mean curvature is bounded  from below by $\mathscr{H}_\phi\ge -l$, $ l \ge 0$, then by the weighted Laplacian comparison theorem \ref{thmW} we have 
$\mathscr{L}_\phi \varrho_{\partial M}(x) \le (m-1)kR +l$ for all $x \in B_R(\partial M)$. Using the properties of the function defined by \eqref{eC} as stated in Lemma \ref{lem45} we derive  $\nabla \psi = (\chi'/R)\nabla \varrho_{\partial M}$ which implies
\begin{align}\label{e4.21}
|\nabla \psi|^2/\psi \le C^2_1/R^2.
\end{align}
We also derive $\Delta \psi =  (\chi'/R)\Delta \varrho_{\partial M} +  (\chi''/R) |\nabla  \varrho_{\partial M}|^2$ which as well implies
\begin{align}\label{e4.22}
\mathscr{L}_\phi \psi & = (\chi'/R) \mathscr{L}_\phi \varrho_{\partial M} +  (\chi''/R) |\nabla  \varrho_{\partial M}|^2 \nonumber\\
& \ge - \frac{C_1}{R}[(m-1)kR+l]\chi^{\frac{1}{2}} - \frac{C_2}{R^2} =  - \frac{1}{R^2}[C_1((m-1)kR+l)R   +C_2].           
\end{align}

\subsection{Proof of Theorem \ref{thm13}}

Since $\mathscr{R}ic_\phi^m\ge -(m-1)k$ and $\mathscr{H}_\phi\ge -l$,  $k,l\ge 0$, we will make use of the weighted Laplacian comparison theorem as in the previous section, and consequently, a cut-off function described in Lemma \ref{lem45} will be applied. First, we define a function $G(x,t)=tF(x,t)$. Again in the consideration of Calabi's trick we assume that the function $\psi(x)G(x,t)$ is smoothly supported in $B_{2R}(\partial M)$, where $\psi$ is the function in Lemma \ref{lem45} and $\psi(x)G(x,t)$ is spatially localized. Suppose also that $\psi G$ attains its maximum at the point $(x_0,t_0) \in B_{2R}(\partial M)\times [0,T]$ for any fixed $t>0$. Without loss of generality, we will assume that $\psi(x_0)G(x_0,t_0)>0$, otherwise the result is obviously trivial.  As before, the proof is divided into two cases, according to whether 
$x_0\not\in \partial M$ or $x_0\in \partial M$.

\noindent
{\bf Case 1:} If $x_0\not\in \partial M$,, we assume $x_0\not\in \text{Cut}(\partial M)$ without loss of generality. By the maximum principle we have at $(x_0,t_0)$
\begin{align}\label{e4.23}
\nabla (\psi G)=0, \ \ \mathscr{L}_\phi(\psi G)\le 0\ \ \ \partial_t(\psi G)\ge 0.
\end{align}
Noting that the first item in \eqref{e4.23} yields $\psi \nabla G = -\nabla \psi G$ which implies $\nabla G= -(\nabla \psi/\psi)G$, while the second item in \eqref{e4.23} yields $\psi \mathscr{L}_\phi G + 2\nabla \psi\nabla G + \mathscr{L}_\phi\psi\cdot G \le 0.$
After combining these two items we obtain 
\begin{align}\label{e4.24}
\psi \mathscr{L}_\phi G - 2(|\nabla \psi|^2/\psi) G + \mathscr{L}_\phi\psi\cdot G \le 0.
\end{align}
Adding the third item in \eqref{e4.23} and \eqref{e4.24}, and referring to the cut-off function estimates \eqref{e4.21} and \eqref{e4.22} we get 
\begin{align}\label{e4.25}
\psi (\mathscr{L}_\phi -\partial_t)G & \le (2|\nabla \psi|^2/\psi -  \mathscr{L}_\phi\psi)G\nonumber\\
& \le \frac{1}{R^2}\Big(2C_1^2+C_2+C_1[(m-1)kR+l]R\Big)G=: \Phi G,
\end{align}
where $\Phi=\Phi(C_1,C_2,n,m,k,R,l)$.
At this juncture, Lemma \ref{lem44} can now be applied by taking the evolution of the quantity $G=tF$:
\begin{align}\label{e4.26}
\psi (\mathscr{L}_\phi -\partial_t)G = & t \psi (\mathscr{L}_\phi -\partial_t)F - G/t \nonumber\\
\ge  & \frac{2(1-\epsilon)}{m\delta^2p}\frac{G^2}{t} +\Big[  \frac{2(1-\epsilon)}{m}\frac{(\delta p+\delta-1)^2}{\delta^2p^2}- 2\Big(\frac{1}{\epsilon}-1\Big)\Big] \frac{|\nabla w|^4}{w^4}t   \nonumber\\
& + \frac{4(1-\epsilon)}{m}\frac{(\delta p+\delta-1)}{\delta^2p^2} \frac{|\nabla w|^2}{w^2}G + \frac{2}{p}\nabla G\nabla \log w \nonumber\\ 
& - 2(m-1)k\frac{|\nabla w|^2}{w^2}t + \Omega^H|\nabla w|^2/w^2t  - H_u G - G/t,
\end{align}
where $\Omega^H:= -[H_u - (\alpha/p^2) uH''(u)]$ and $H_u:=H'(u)-H(u)/u$.

\noindent
By the Cauchy-Schwarz  and Young's  inequalities we estimates
\begin{align}\label{e4.27}
2(m-1)k\frac{|\nabla w|^2}{w^2}t \le \frac{2(1-\epsilon)(\delta-1)^2}{m\delta^2p^2}\frac{|\nabla w|^4}{w^4}t + \frac{m\delta^2p^2k^2}{2(1-\epsilon)(\delta-1)^2}t,
\end{align}
\begin{align}\label{e4.28}
 \Omega_H|\nabla w|^2/w^2t \ge - \frac{(\delta-1)^2}{2}\frac{|\nabla w|^4}{w^4}t - \frac{1}{(\delta-1)^2}\sup_{\mathcal{Q}_{2R,T}}  \{[\Omega^H]_+\}^2 t.
\end{align} 
In the next we denote by $\mu^H$ and $\gamma^H$ the constants  
$$\mu^H =   \sup_{\mathcal{Q}_{2R,T} } \{[H_u]_+\}  \ \ \text{and}\ \ \gamma^H= \sup_{\mathcal{Q}_{2R,T} }  \{[\Omega^H]_+\}.$$ 
Choosing $\delta>1$ and $p>0$ such that 
$1/p >\Big[ (m\delta)/[(1-\epsilon)(\delta -1)]\Big]\Big[1/\epsilon-1 + (\delta-1)^2/4\Big]$ satisfies the hypothesis \eqref{hyp5} and again makes
\begin{align}\label{e4.29}
\Big[  \frac{2(1-\epsilon)}{m}\frac{(\delta-1)}{\delta p}- 2\Big(\frac{1}{\epsilon}-1\Big)- \frac{(\delta-1)^2}{2}\Big] \frac{|\nabla w|^4}{w^4}t  \ge 0.
\end{align}
Now plugging in \eqref{e4.27} - \eqref{e4.29} into \eqref{e4.26}, noting some cancellation and multiplying through by $\psi$ gives
\begin{align}\label{e4.30}
\psi (\mathscr{L}_\phi -\partial_t)G \ge  & \frac{2(1-\epsilon)}{m\delta^2p}\frac{G^2}{t}\psi  + \frac{4(1-\epsilon)}{m}\frac{(\delta p+\delta-1)}{\delta^2p^2} \frac{|\nabla w|^2}{w^2}G \psi  \nonumber\\
& + \frac{2}{p}\psi \nabla G\nabla \log w  -  \frac{m\delta^2p^2k^2}{2(1-\epsilon)(\delta-1)^2}t \nonumber \\ 
& -  \frac{1}{(\delta-1)^2}{\gamma^H}^2\psi t - \mu^H G\psi - \psi G/t.
\end{align}
By Young's inequality again while using the identity $\psi \nabla G=-G\nabla \psi$ derived from \eqref{e4.23} we have 
\begin{align}\label{e4.31}
\frac{2}{p}\psi \nabla G\nabla \log w  & = - \frac{2}{p}G \nabla\psi \frac{\nabla w}{w}\nonumber \\
& \ge - \frac{4(1-\epsilon)}{m}\frac{(\delta p+\delta-1)}{\delta^2p^2} \frac{|\nabla w|^2}{w^2}G \psi   - \frac{m\delta^2}{4(1-\epsilon)(\delta p+\delta-1)} \frac{|\nabla \psi|^2}{\psi}G.
\end{align}
Recall that the quantity appearing on the left hand side of \eqref{e4.30} has been estimated in \eqref{e4.25}. Therefore combining \eqref{e4.25},  \eqref{e4.30} and \eqref{e4.31}  and rearranging terms leads to
\begin{align*}
 \frac{2(1-\epsilon)}{m\delta^2p}\frac{G^2}{t}\psi   - \frac{m\delta^2}{4(1-\epsilon)(\delta p+\delta-1)}& \frac{|\nabla \psi|^2}{\psi}G - \frac{m\delta^2p^2k^2}{2(1-\epsilon)(\delta-1)^2}t \nonumber \\
 &   - \frac{1}{(\delta-1)^2}{\gamma^H}^2 \psi t - \mu^H  G\psi - \psi G/t  - \Phi G \le 0.
\end{align*}
To obtain a more compact expression from the last inequality, multiply through it by $\psi t$ again and rearrange terms we obtain
\begin{align}\label{e4.32}
 \frac{2(1-\epsilon)}{m\delta^2p} G^2 \psi^2 -(1+ \mu^H t +\widetilde{\Phi}t)G\psi 
 & - t^2 \Big[\frac{m\delta^2p^2k^2}{2(1-\epsilon)(\delta-1)^2} + \frac{1}{(\delta-1)^2}{\gamma^H}^2 \Big]\psi^2 \le 0,
\end{align}
where
$\widetilde{\Phi} = \frac{C_1^2}{R^2}\frac{m\delta^2}{4(1-\epsilon)(\delta p+\delta-1)} + \Phi$ with $\Phi$ being as defined in \eqref{e4.25}.

Observe that the last inequality is obviously quadratic in $\psi G$ and can therefore be compared with the following quadratic inequality $ay^2-by-c\le 0$ with $a>0$, $b,c\ge 0$ being real coefficients, which can be solved using the formula $y=b/a+\sqrt{c/a}$. So setting $y=\psi G$, $a= 2(1-\epsilon)/(m\delta^2p^2)$, $b= (1+ \mu^H t +\widetilde{\Phi}t)$ and 
$$c= \frac{m\delta^2p^2k^2}{2(1-\epsilon)(\delta-1)^2} + \frac{1}{(\delta-1)^2}{\gamma^H}^2 $$
we solve \eqref{e4.32} for $\psi G$ at $(x_0,t_0)$ and obtain 
\begin{align*}
\psi G & \le \frac{m\delta^2p^2}{2(1-\epsilon)} (1+ \mu^H t_0+\widetilde{\Phi}t_0) + t_0 \sqrt{ \frac{m\delta^2p^2}{2(1-\epsilon)}}\sqrt{\Big(\frac{m\delta^2p^2k^2}{2(1-\epsilon)(\delta-1)^2} +  \frac{1}{(\delta-1)^2}{\gamma^H}^2 \Big)\psi^2}\\
&\le \frac{m\delta^2p^2}{2(1-\epsilon)} (1+ \mu^H T+\widetilde{\Phi}T) + T \sqrt{ \frac{m\delta^2p^2}{2(1-\epsilon)}} \left(\sqrt{\frac{m\delta^2p^2k^2}{2(1-\epsilon)(\delta-1)^2}} + \frac{1}{(\delta-1)}\gamma^H \right)\psi.
\end{align*}
Taking into consideration the properties of $\psi$ (see Lemma \ref{lem45}) we set $\psi\equiv 1$ in $B_{2R}$ and we finally write
\begin{align*}
& \sup G(x,t)  \le \sup \{\psi(x_0)G(x_0,t_0)\}\\
&\le  \frac{m\delta^2p^2}{2(1-\epsilon)} (1+ \mu^H T+\widetilde{\Phi}T) + T \sqrt{ \frac{m\delta^2p^2}{2(1-\epsilon)}} \left(\sqrt{\frac{m\delta^2p^2k^2}{2(1-\epsilon)(\delta-1)^2}} + \frac{1}{(\delta-1)}\gamma^H \right)\psi.
\end{align*}
at $(x_0,t_0)\in \mathcal{Q}_{2R,T}(\partial M)$.

Noting that $G=tF$ with $F$ defined by $F= |\nabla w|^2/w^2 +\alpha \widehat{H}(w)+\beta w_t/w$, $\alpha,\beta>0$ and by choice $\beta=\alpha/p$ and $\alpha=\delta p^2$.  Thus, referring to $w=u^{-p}$ and $\widehat{H}(w)=H(u)/u$ we can write 
$$G/t = F = p^2|\nabla u|^2/u^2 + \delta p^2 H(u)/u - \delta p^2 u_t/u.$$
Therefore  
\begin{align*}
\frac{|\nabla u|^2}{u^2} &+\delta\frac{H(u)}{u} -\delta \frac{u_t}{u}\\
& \le \frac{m\delta^2}{2(1-\epsilon)}\Big(\frac{1}{T}+ \mu^H+ \widetilde{\Phi}\Big) + \frac{m\delta^2 k}{2(1-\epsilon)(1-\delta)} +  \frac{1}{p^2(\delta-1)}\sqrt{ \frac{m \delta^2 p^2}{2(1-\epsilon)}}\gamma^H. 
\end{align*}
at $(x,t)\in \mathcal{Q}_{R,T}(\partial M)$. Therefore \eqref{eqthm3} follows since $T>0$ was arbitrarily fixed and the theorem is proved for the case $x\not\in \partial M$.

\noindent
{\bf Case 2:}
If $x_0 \in \partial M$, we will show that the gradient estimates \eqref{eqthm3} holds on $\partial M \times [0, T]$. As in the discussion ensued in  previous sections, at the maximum point $(x_0,t_0)$ over  $\partial M \times [0, T]$ we have
\begin{align}\label{f1}
F_\nu  = (|\nabla w|^2/w^2)_\nu \ge 0
\end{align}
since $w_t=0$ and $\alpha \widehat{H} +\beta w_t/w \ge 0$ on $\partial M \times [0, T]$ as informed by the theorem assumption $u_t \le H(u)$,  $u_t=0$ on $\partial M \times [0, T]$. This can be seen by setting  $w=u^{-p}$ with $\beta = \alpha/p\ge 0$.

Noting also that $|\nabla w| = w_v \ge 0$ since $u>0$ and $u_v\ge0$. By the Reilly-type formula \eqref{Rei} we have
\begin{align}\label{f2}
(|\nabla w|^2/w^2)_\nu  = 2(w_\nu/w^2)[\mathscr{L}_\phi w-\mathscr{H}_\phi w_\nu] - 2w_\nu^2w_\nu/w^3.
\end{align}
From Lemma \ref{lem41} we have $\mathscr{L}_\phi w = ((p+1)/p)(w_\nu^2/w) + w_t +pw\widehat{H}$. It then follows by combining \eqref{f1} and \eqref{f2} that 
\begin{align*}
0\le (w_\nu/w^2)[(1/p)(w_\nu^2/w) + w_t +pw\widehat{H}-\mathscr{H}_\phi w_\nu]
\end{align*}
at $(x_0,t_0)$. Obviously the condition $u_t\le H(u)$ implies $w_t+pw\widehat{H}\ge 0$ over  $\partial M \times [0, T]$. Hence
\begin{align*}
-(1/p)(w_\nu/w)  + \mathscr{H}_\phi \le 0.
\end{align*}
Using the fact that $w_\nu/w = -pu_\nu/u$ and $p>0$ we obtain
$$u_v+ \mathscr{H}_\phi \le 0$$
at $(x_0,t_0)$.  This implies that 
$$F(x_0,t_0)\le l^2$$
which holds for all $(x,t) \in \mathcal{Q}_{R,T}(\partial M)$ by noting that $\sup F(x,t)\le \sup F(x_0,t_0)$.
This therefore completes the proof.

\qed

\section{Applications of Li-Yau  type gradient estimates (Theorem \ref{thm13})}\label{sec6}
\subsection{Global Li-Yau gradient estimates}
Theorem \ref{thm13} which is local in nature implies the following global gradient estimates under appropriate global bounds as seen in the following theorem.
\begin{theorem}\label{thm61}
Let $M$ be an $n$-dimensional complete smooth metric measure space with compact boundary whose mean curvature satisfies the bound $ \mathscr{H}_\phi(\partial M) \ge - l$, $l \ge 0$. Let $u=u(x,t)$ be a smooth positive solution to \eqref{11} such that $u$ satisfies the Dirichlet boundary condition (i.e., $u(\cdot,t)|_{\partial M}$ is constant for each time $t$), 
$u_\nu\ge 0$ and $u_t \le H(u)$ on $\partial M \times (-\infty,+\infty)$.  Denote
$$\mu^H : =   \sup_{M\times (-\infty,+\infty)} \{[H_u]_+\}  \ \ \text{and}\ \ \gamma^H := \sup_{M\times (-\infty,+\infty)}  \{[\Omega^H]_+\},$$ 
where the quantities $[H_u]_+$ and $[\Omega^H]_+ $ are as defined in Theorem \ref{thm13}.
\begin{enumerate}
\item Suppose $\mathscr{R}ic_\phi^m \ge -(m-1)k, \ k\ge 0$, then for all $x\in M$ and $t\neq 0$, we have
\begin{align}\label{eq61}
\frac{|\nabla u|^2}{u^2} &+\delta\frac{H(u)}{u} -\delta \frac{u_t}{u}
 \le \frac{m\delta^2}{2}\Big(\frac{1}{t}+ \frac{k}{\delta-1}+\mu^H\Big) +  \frac{\delta}{p(\delta-1)} \sqrt{ \frac{m}{2}}\gamma^H,
\end{align}
where $\delta>1$ and  $p>0$.
\item Suppose $\mathscr{R}ic_\phi^m \ge  0$, then for all $x\in M$ and $t\neq 0$, we have
\begin{align}\label{eq62}
\frac{|\nabla u|^2}{u^2} &+ \delta \frac{H(u)}{u} -  \delta \frac{u_t}{u}
 \le \frac{m\delta^2}{2t} + \frac{m\delta^2}{2}\mu^H +  \frac{\delta}{p(\delta-1)} \frac{ \sqrt{m}}{ \sqrt{2}}\gamma^H,
\end{align}
where $\delta>1$ and $p>0$.
\end{enumerate}
\end{theorem}

\proof
The proof is a direct consequence of Theorem \ref{thm13} once the bounds on $\mathscr{R}ic_\phi^m$ and $\mathscr{H}_\phi$, and all the 
supremums are  taken in global sense.  Estimates \eqref{eq61} therefore follows from \eqref{eqthm3} by sending $\epsilon  \searrow 0$, whilts \eqref{eq62} follows also by setting  $k\equiv 0$  since $\mathscr{R}ic_\phi^m \ge  0$ for all $x\in M$.
\qed

\subsection{Li-Yau Harnack inequalities}
The gradient estimates obtained in Theorem \ref{thm13} can be utilized to derive parabolic Harnack inequality for the positive solution of \eqref{11} on $M\times (-\infty,\infty)$ as seen in the following theorem.
\begin{theorem}\label{thm62}
Let $M$ be an $n$-dimensional complete smooth metric measure space with compact boundary  and let $u$ be a smooth positive solution to \eqref{11}. Under the assumptions of Theorem \ref{thm13} we have for all $(x,t) \in \mathcal{Q}_{R,T}(\partial M)$ 
with $x_1,x_2 \in B_R(\partial M)$, $-\infty <t_1<t_2< +\infty$, $t_1, t_2\neq 0$ and $\delta>1$ that 
\begin{align}\label{eq63}
u(x_1,t_1)\le u(x_2,t_2)\left(\frac{t_2}{t_1}\right)^{\frac{m\delta}{2(1-\epsilon)}}\exp \left\{\frac{\delta d^2(x_1,x_2)}{4(t_2-t_1)} + A (t_2-t_1) \right\},
\end{align}
where $d(x_1,x_2)$ is the geodesic distance between $x_1$ and $x_2$, and 
\begin{align*}
A &=A(m,p,\delta,\epsilon,k,\mu^H,\gamma^H,\sigma^H) \\
& =   \frac{m\delta}{2(1-\epsilon)}\Big(\frac{k}{\delta-1}+ \widetilde{\Phi}  +\mu^H\Big) +  \frac{1}{p(\delta-1)} \sqrt{ \frac{m}{2(1-\epsilon)}}\gamma^H - \sigma^H
\end{align*}
with $\displaystyle \sigma^H = \inf_{\mathcal{Q}_{2R,T}(\partial M)} \left\{[H(u)/u]_-\right\}$, $\delta>1$ and $p>0$.
\end{theorem}

\proof
For the purpose of the proof of Li-Yau Harnack inequality, \eqref{eqthm3} can be written as 
\begin{align}\label{eq64}
- \frac{u_t}{u} \le \frac{m\delta}{2(1-\epsilon)t} - \frac{|\nabla u|^2}{\delta u^2} + A,
\end{align}
where
\begin{align*}
A& =  \frac{m\delta}{2(1-\epsilon)}\Big(\frac{k}{\delta-1}+ \widetilde{\Phi}  +\mu^H\Big) +  \frac{1}{p(\delta-1)} \sqrt{ \frac{m}{2(1-\epsilon)}}\gamma^H - \sigma^H,\\
 \sigma^H & = \inf_{\mathcal{Q}_{2R,T}(\partial M)} \left\{[H(u)/u]_-\right\}.
 \end{align*}
Let $\gamma$ be the shortest geodesic joining $x_1$ and $x_2$ in $\mathcal{Q}_{R,T}(\partial M)$ such that $\gamma(t_1)=x_1$ and $\gamma(t_2)=x_2$. By a straightforward computation we have 
\begin{align*}
\log u(\gamma(t),t)\Big|_{t_1}^{t_2} & = \int_{t_1}^{t_2} \frac{d}{dt}\log u(\gamma(t),t) dt\\
& = \int_{t_1}^{t_2} \left[\dot{\gamma} \nabla u/u + u_t/u \right]dt\\
& \ge \int_{t_1}^{t_2} \left[\dot{\gamma} \nabla u/u + |\nabla u|^2/(\delta u^2) -A -m\delta /2(1-\epsilon)t \right]dt
\end{align*}
by using \eqref{eq64}.   Applying the completing the square method to rewrite the expression $\dot{\gamma} \nabla u/u + |\nabla u|^2/(\delta u^2) $ and using the fact that $(1/\delta)\left| (\delta/2)\dot{\gamma} + \nabla u/u \right|^2\ge 0$,  we have
\begin{align*}
\log u(x_2,t_2)-\log u(x_1,t_1) \ge  \int_{t_1}^{t_2} \left[-\frac{\delta}{4}|\dot{\gamma}|^2 - A -  \frac{m\delta}{2(1-\epsilon)t}\right]dt.
\end{align*}
Therefore
\begin{align}\label{eq65}
\log \frac{u(x_2,t_2)}{u(x_1,t_1)} \ge -  \int_{t_1}^{t_2}\frac{\delta}{4}|\dot{\gamma}|^2 dt - A(t_2-t_1) -   \frac{m\delta}{2(1-\epsilon)}\log\left(\frac{t_2}{t_1}\right).
\end{align}
By exponentiating both sides of \eqref{eq65} we arrive at
\begin{align*}
u(x_2,t_2)/u(x_1,t_1) \ge \left(t_2/t_1\right)^{-\frac{m\delta}{2(1-\epsilon)}}\exp \left\{\delta/4  \int_{t_1}^{t_2} |\dot{\gamma}|^2 dt   -A (t_2-t_1) \right\},
\end{align*}
which when rearranging yields
\begin{align*}
u(x_1,t_1) \le u(x_2,t_2) \left(t_2/t_1\right)^{\frac{m\delta}{2(1-\epsilon)}}\exp \left\{\delta \rho + A (t_2-t_1) \right\},
\end{align*}
where $\displaystyle \rho = \rho(x_1,x_2, t_2-t_1) = \inf_\gamma \left\{\frac{1}{4(t_2-t_2)} \int_0^1 |\dot{\gamma}|^2 dt \right\}.$
We arrive at the desired inequality by using the fact that $d(x_1,x_2) =\inf_{\gamma\in \Gamma}  \int_0^1 |\dot{\gamma}| dt$ since every pair of points in $M$ can be joined by the minimizing geodesic, with
$$\Gamma:=\left\{\gamma:[0,1]\to M \ | \ \gamma(0)=x_1, \gamma(1)=x_2\right\}$$
being a smooth curve in $\mathcal{Q}_{2R,T}(\partial M)$.

\qed

Theorem \ref{thm62} gives the local Li-Yau inequality, which also implies some global form if the bound are taken in global sense with the same argument. For instance we obtain the following corollary.
\begin{corollary}\label{cor63}
Let $M$ be an $n$-dimensional complete smooth metric measure space with compact boundary satisfying $\mathscr{R}ic_\phi^m \ge 0$ and $\mathscr{H}_\phi \ge 0$.  Let $u$ be a smooth positive solution to \eqref{11} satisfying the assumption of Theorem \ref{thm61}. Then we have for all $x_1,x_2\in M$ $-\infty <t_1<t_2< +\infty$, $t_1, t_2\neq 0$ that 
\begin{align}\label{eq66}
u(x_1,t_1)\le u(x_2,t_2)\left(\frac{t_2}{t_1}\right)^{\frac{m}{2}}\exp \left\{\frac{\delta d^2(x_1,x_2)}{4(t_2-t_1)} + \widetilde{A} (t_2-t_1) \right\},
\end{align}
$$A=A(m,p,\delta, \mu^H,\gamma^H,\sigma^H) =   \frac{m}{2}\mu^H + \frac{1}{p(\delta-1)}\frac{\sqrt{m}}{\sqrt{2}}\gamma^H  - \sigma^H$$
with $\displaystyle \sigma^H = \inf_{M\times (-\infty,+\infty)} \left\{[H(u)/u]_-\right\}$, $\delta>1$ and $p>0$.
\end{corollary}

\proof
The proof of Corollary \ref{cor63} follows the same argument as in the proof of Theorem \ref{thm62} with the use of global bounds and global gradient estimates in Theorem \ref{thm61}. Since 
$\mathscr{R}ic_\phi^m \ge 0$ we just set $k\equiv 0$ and  $\epsilon \searrow 0$, and we obtain \eqref{eq66} at once.

\qed

\begin{remark}
If we set $\epsilon = 1/2$ we still obtain \eqref{eq66} with power of $(t_2/t_1)$ becoming $m$, while the quantity $A$ becomes
$A=m\mu^H+ (\sqrt{m}/p)\gamma^H -\sigma^H$.
\end{remark}

\subsection{Li-Yau type gradient estimates for elliptic equations and Liouville-type theorems}
Here we discuss Li-Yau type gradient estimates for elliptic form (steady state) of \eqref{11}. The first result here is the local gradient estimates on the positive smooth solution of the following weighted elliptic equation
\begin{align}\label{eq67}
\begin{cases}
\displaystyle \mathscr{L}_\phi u +H(u) &= 0, \\
\displaystyle \hspace{1.5cm} u|_{\partial M} & = \text{const.},
\end{cases}
\end{align}
where $H=H(s)$ is a sufficiently smooth nonlinear function for $s>0$.

\begin{theorem}\label{thm67}
Let $M$ be an $n$-dimensional complete smooth metric measure space with compact boundary satisfying $\mathscr{R}ic_\phi^m \ge -(m-1)k$ and $\mathscr{H}_\phi \ge -l$, $k,l \ge 0$.  Let $u=u(x)$ be a smooth positive solution to \eqref{eq67} on $B_R(\partial M)$ such that $u_\nu\ge 0$ over $\partial M$.  Then for some $\epsilon \in (0,1)$ and $\delta>1$, we have
\begin{align}\label{eq68}
\frac{|\nabla u|^2}{u^2} &+\delta\frac{H(u)}{u}
 \le \frac{m\delta^2}{2(1-\epsilon)}\Big(\frac{k}{\delta-1}+\widetilde{\Phi}+  \widetilde{\mu}^H\Big) +  \frac{\delta}{p(\delta-1)} \sqrt{ \frac{m}{2(1-\epsilon)}}\widetilde{\gamma}^H,
\end{align} 
for all $x\in B_R(\partial M)$ with 
$0<p \le [(1-\epsilon)(\delta-1)]/\left[m\delta \left(\frac{1}{\epsilon} -1 + \frac{(\delta-1)^2}{4} \right)\right]$, where
$$\widetilde{\mu}^H : =   \sup_{B_R(\partial M)} \{[H_u]_+\}  \ \ \text{and}\ \ \widetilde{\gamma}^H := \sup_{B_R(\partial M)}  \{[\Omega^H]_+\},$$ 
$[H_u]_+$ and $[\Omega^H]_+$ are as defined in Theorem \ref{thm13}.
\end{theorem}

Globally we have the following Li-Yau type gradient estimates.
\begin{theorem}\label{thm66}
Let $M$ be an $n$-dimensional complete smooth metric measure space with compact boundary satisfying $\mathscr{R}ic_\phi^m \ge0$ and $\mathscr{H}_\phi \ge 0$.   Let $u=u(x)$ be a smooth positive solution to \eqref{eq67} such that $u_\nu\ge 0$ over $\partial M$.  Then for some $\delta>1$, we have
\begin{align}\label{eq69}
\frac{|\nabla u|^2}{u^2} &+\delta\frac{H(u)}{u}
 \le \frac{m\delta^2}{2(1-\epsilon)} \widetilde{\mu}_M^H +  \frac{\delta}{p(\delta-1)} \sqrt{ \frac{m}{2(1-\epsilon)}}\widetilde{\gamma}_M^H,
\end{align} 
for all $x\in M$, where
$$\widetilde{\mu}_M^H : =   \sup_M \{[H_u]_+\}  \ \ \text{and}\ \ \widetilde{\gamma}_M^H := \sup_M  \{[\Omega^H]_+\}.$$ 
In particular, choosing 
$p=[(1-\epsilon)(\delta-1)]/\left[m\delta \left(\frac{1}{\epsilon} -1 + \frac{(\delta-1)^2}{4} \right)\right]$,  $\delta=2$ and $\epsilon=\frac{1}{2}$, \begin{align}\label{eq70}
\frac{|\nabla u|^2}{u^2} &+ 2\frac{H(u)}{u}
 \le 4m \ \widetilde{\mu}_M^H +  10m\sqrt{m}\ \widetilde{\gamma}_M^H.
\end{align} 
\end{theorem}

As a consequence of the above  Li-Yau type gradient estimates we can prove the following Liouville-type results on nonnegative weighted Ricci curvature and mean curvature conditions.

\begin{theorem}\label{thm67}
Let $M$ be an $n$-dimensional complete smooth metric measure space with compact boundary satisfying $\mathscr{R}ic_\phi^m \ge0$ and $\mathscr{H}_\phi \ge 0$.   Let $u=u(x)$ be a smooth bounded positive solution to \eqref{eq67} such that $u_\nu\ge 0$ over $\partial M$, where $H(u)\ge 0$, $H'(u)-H(u)/u\le 0$ and $H'(u)-H(u)/u -\delta uH''(u)\ge 0$ for $\delta>1$ and all $x \in M$. Then $u$ must be a constant.
\end{theorem}

\subsection{Some special examples of $H(u)$}\label{sec64}
Lastly in this subsection, we demonstrate that the global estimates in Theorem \ref{thm66} can be applied to investigate Liouville-type results for more special cases of $H(u)$ with  a power like growth or logarithmic-type or exponential-type nonlinearities, such as cases of $H(u)= au^\lambda (\log u)^{\alpha}$,  $H(u)= au^q + bu^s$, $H(u)= Ae^{au}$ or  a mixture of these nonlinearities.

\subsection*{1. The case $H(u)= au^\lambda \log u$ with $ a<0$, $\lambda\ge 1$}
Consider $H(u)=au^\lambda (\log u)^{\alpha}$, but for simplicity we keep $\alpha=1$, $a\le 0$ and $\lambda\ge 1$. Thus we consider the specific case $H(u)=au^\lambda \log u$ for bounded positive function $u$. Here we see that 
$$H_u=H'(u)-H(u)/u = a\left(1+(\lambda-1)\log u\right)u^{\lambda-1}$$
and 
$$H'(u)-H(u)/u-\delta uH''(u) = a\left([1-\delta(2\lambda-1)]+ [(\lambda-1)(1-\delta\lambda)]\log u \right)u^{\lambda-1}$$
which satisfy the assumptions that $H'(u)-H(u)/u \le 0$ and $H'(u)-H(u)/u-\delta uH''(u)\ge 0$ for $0<u(x)<c\le 1$ since $\delta >1$, $a\le 0$ and $\lambda\ge 1$. Clearly
$$\widetilde{\mu}^H_M = \sup_M \{[H_u]_+\} =[a(1+(\lambda-1)\log u)]_+ \sup_M\{u^{\lambda-1}\} =0$$
and
$$\widetilde{\gamma}_M^H = \sup_M  \{[\Omega^H]_+\} =  \sup_M  \{[-[H'(u)-H(u)/u-\delta uH''(u)]]_+\} =0.$$ 
By the above discussion we have the following Liouville-type results
\begin{proposition}
Let $M$ be an $n$-dimensional complete smooth metric measure space with compact boundary satisfying $\mathscr{R}ic_\phi^m \ge0$ and $\mathscr{H}_\phi \ge 0$.   Let $u$ be a bounded positive solution to 
\begin{align}
\begin{cases}
\displaystyle \mathscr{L}_\phi u + au^\lambda \log u & =0, \hspace{.5cm} \lambda \ge 1\\
\displaystyle \hspace{2cm} u |_{\partial M} & = \text{const.}
\end{cases}
\end{align}
such that $u_\nu\ge 0$ over $\partial M$. If $a\le 0$, then $u$ must be a constant. If a is further restricted to $a<0$, then $u \ge \exp\{2/(1-\lambda)\}$.
\end{proposition}

\proof
Combining the above discussion (which gives $\widetilde{\mu}^H_M =0$,  $\widetilde{\gamma}_M^H =0$) together with \eqref{eq70} (Theorem \ref{thm66}), it then follows that $|\nabla u|\equiv 0$. Referring to the assumption that $H(u)\ge 0$ and $H'(u)-H(u)/u-\delta uH''(u) \ge 0$, we have $\log u\ge 2/(1-\lambda)$.

\qed

\subsection*{2. The case $H(u)= au^q + bu^s$ with real coefficients $a,b$ and real exponents $s,q$}
Consider $H(u)=au+bu^q$ for bounded positive function $u$. Here we compute that 
$$H_u=H'(u)-H(u)/u = b(q-1)u^{q-1}$$
and 
$$H'(u)-H(u)/u-\delta uH''(u) = b(q-1)(1-\delta q)u^{q-1}$$
which satisfy the assumptions that $H'(u)-H(u)/u \le 0$ and $H'(u)-H(u)/u-\delta uH''(u)\ge 0$ for $\delta >1$ and for either $b\ge 0$ and $q<1$ or $b\le 0$ and $q>1$. Clearly
$$\widetilde{\mu}^H_M = \sup_M \{[H_u]_+\} =[b(q-1)]_+ \sup_M\{u^{q-1}\} =0$$
and
$$\widetilde{\gamma}_M^H = \sup_M  \{[\Omega^H]_+\} = [ b(q-1)(\delta q-1)]_+ \sup_M  \{u^{q-1}\} =0.$$ 
The following Liouville-type results can now be proved in a similar way.
\begin{proposition}
Let $M$ be an $n$-dimensional complete smooth metric measure space with compact boundary satisfying $\mathscr{R}ic_\phi^m \ge0$ and $\mathscr{H}_\phi \ge 0$.   Let $u$ be a bounded positive solution to 
\begin{align}\label{eqac}
\begin{cases}
\displaystyle \mathscr{L}_\phi u + au + bu^q & =0, \\
\displaystyle \hspace{2cm} u |_{\partial M} & = \text{const}.
\end{cases}
\end{align}
on $M$ such that $u_\nu\ge 0$ over $\partial M$, where $a,b,q\in \mathbb{R}$ and $q\neq 1$.  If either $b\ge 0$ and $q<1$ or $b\le 0$ and $q>1$,  then $u$ must be a constant.  Moreover, if $a>0$, $b<0$ and $q>1$, then $u$ is the constant $u=\sqrt[q-1]{a/(-b)}$.  
\end{proposition}
\begin{remark}
The last proposition includes Liouville-type results for bounded solution of the following specific equation of mathematical physics, namely;
\begin{enumerate}
\item Weighted Allen-Cahn equation: $\mathscr{L}_\phi u + u(1-u^3) =0$, $0<u\le 1$. That is, setting $a=1$, $b=-1$ and $q=3$ in \eqref{eqac}.
\item Weighted elliptic Fisher KPP equation: $\mathscr{L}_\phi u + cu(1-u) =0$, $c>0$, $0<u\le 1$. That is, setting $a=c>0$, $b=-c$ and $q=2$ in \eqref{eqac}.
\item Weighted semilinear elliptic equation: $\mathscr{L}_\phi u +bu^q=0$, $b, q$ are real. That is, setting $a=0$ in \eqref{eqac}.
\item In general, a similar Liouville type result holds for the case $\displaystyle H(u)=\sum_{l=1}^r b_lu^{\lambda_l}$, $1\le l\le r$. That is, any bounded positive solution to the weighted boundary value problem
\begin{align}\label{eqac}
\begin{cases}
\displaystyle \mathscr{L}_\phi u + \sum_{l=1}^r b_lu^{\lambda_l} & =0, \\
\displaystyle \hspace{2cm} u |_{\partial M} & = \text{const}.
\end{cases}
\end{align}
on $M$ such that $u_\nu\ge 0$ over $\partial M$ with $\mathscr{R}ic_\phi^m\ge 0$ and $\mathscr{H}_\phi\ge 0$ must be constant if either $b_l's\ge 0$ and $\lambda_l\le 1$ or $b_l's\le 0$ and $\lambda_l\ge 1$.
\end{enumerate}
\end{remark}

\section{Applications of elliptic gradient estimates}\label{sec7}
First in this section we mention that local Hamilton-Souplet-Zhang type gradient estimates in Theorem \ref{thm11} (resp. Theorem \ref{thm21a}) and Theorem \ref{thm21} (resp. Theorem \ref{thm22}) have global counterparts when ensuring that all bounds and supremums are globally defined on the whole of $M$. The resulting global gradient estimates can be proved by allowing $R$ to pass to the limit $R\to \infty$ in the respective local estimates. To this end and for the purpose of other applications of elliptic gradient estimates, we focus on Theorem \ref{thm21}, while noting that Theorem \ref{thm11} can be applied in similar contexts.

\subsection{Global elliptic gradient estimates}
Theorem \ref{thm21} yields the following global estimates by passing to the limit $R\to \infty$.

\begin{theorem}\label{thm71}
Let $M$ be an $n$-dimensional complete smooth metric measure space with compact boundary satisfying $\mathscr{R}ic_\phi^m \ge -(m-1)k$ and $ \mathscr{H}_\phi \ge - l$ for some $k, l \ge 0$. Let $u$ be a smooth positive solution to \eqref{11} such that $u$ satisfies the Dirichlet boundary condition,
$u_\nu\ge 0$ and $u_t \le H(u)$ on $\partial M \times (-\infty,+\infty)$.  
Then there exists a constant $C_\varepsilon>0$ depending only on $n,m, \varepsilon$ with  $\varepsilon \in (0,  1/2)$ such that 
\begin{align}\label{eq71}
\frac{|\nabla u|}{u^{1-3\varepsilon/2}} & \le C_\varepsilon\Bigg\{ \frac{1}{\sqrt{T}}+\sqrt{L}+ \Gamma^\varepsilon_+[u] \Bigg\}\Big(\sup_{M\times (-\infty,\infty)}\{\sqrt{u^{3\varepsilon}}\}\Big) ,
\end{align}
where  $L$ and $\Gamma^\varepsilon_+[u]$ are as defined in Theorem \ref{thm21}, but with the supremums taking all over the whole of $M$.  
\end{theorem}

An immediate consequence of the above result is the following global gradient estimates for weighted nonlinear elliptic equation on manifold with boundary
\begin{align}\label{eq72}
\displaystyle \mathscr{L}_\phi u(x) +H(u(x)) &= 0 \ \text{in}\ M,   \ \ \ u(x)  = \text{const.} \ \text{on} \ \  \partial M.
\end{align}
\begin{theorem}\label{thm72}
Let $M$ be an $n$-dimensional complete smooth metric measure space with compact boundary satisfying $\mathscr{R}ic_\phi^m \ge -(m-1)k$ and $ \mathscr{H}_\phi \ge - l$, $k, l \ge 0$. Let $u=u(x)$ be a smooth positive bounded solution to \eqref{eq72} such that
$u_\nu\ge 0$ over $\partial M$.  Then for $\varepsilon \in (0,1/2]$ we have 
\begin{align}\label{eq73}
\frac{|\nabla u|}{u^{1-3\varepsilon/2}} & \le C_\varepsilon\Big(\sqrt{L}+ \Gamma^\varepsilon_+[u] \Big)\Big(\sup_{M}\{\sqrt{u^{3\varepsilon}}\}\Big) ,
\end{align}
where the constant $C_\varepsilon$ and the quantities $L$ and $\Gamma^\varepsilon_+[u]$ are as defined in Theorem \ref{thm71}.
\end{theorem}

\proof
The proof is an immediate consequence of Theorem \ref{thm71} by letting $T\to \infty$ since the solution of \eqref{eq72} is time-independent.
\qed

\subsection{Liouville-type results}
Another important implication of these estimates is the following Liouville-type theorem.
 \begin{theorem}\label{thm73}
Let $M$ be an $n$-dimensional complete smooth metric measure space with compact boundary satisfying $\mathscr{R}ic_\phi^m \ge0$ and $\mathscr{H}_\phi \ge 0$.   Let $u=u(x)$ be a smooth positive bounded solution to \eqref{eq72} such that $u_\nu\ge 0$ over $\partial M$, where  $2H'(u)+(3\varepsilon-2)H(u)/u\le 0$ for all $x \in M$ and $\varepsilon \in (0,1/2]$.
Then $u$ must be constant.  In particular $H(u)\equiv 0$.
 \end{theorem}

\proof
The proof follows from Theorem \ref{thm72}. By the assumption $2H'(u)+(3\varepsilon-2)H(u)/u\le 0$,  $\varepsilon \in (0,1/2]$, we see that 
$$ \Gamma^\varepsilon_+[u] = \left[2H'(u)+(3\varepsilon-2)H(u)/u \right]_+=0.$$
Considering also the global bounds $\mathscr{R}ic_\phi^m \ge0$ and $\mathscr{H}_\phi \ge 0$ we set $L=(k^2+l^4)^{1/2}=0$. Thus the global estimates in Theorem \ref{thm72} implies
\begin{align*}
\frac{|\nabla u|}{u^{1-3\varepsilon/2}} & \le C_\varepsilon \Gamma^\varepsilon_+[u]\Big(\sup_{M}\{\sqrt{u^{3\varepsilon}}\}\Big) = 0
\end{align*}
since $\sup_M \{\sqrt{u^{3\varepsilon}}\}$ is finite. It therefore follows that $|\nabla u|\equiv 0$ on $M$, and by implication $u$ is a constant. Substituting this constant solution into \eqref{eq71} yields $H(u)\equiv 0$.

\qed

Finally, we remark that the above global gradient estimates can be applied to establish Liouville-type theorems for several special cases of $H(u)$ as in the previous section. As examples we want to consider two specific nonlinear elliptic equations,  namely:
\begin{align}\label{eq74}
\displaystyle \mathscr{L}_\phi u + \sum_{j=1}^r q_ju^{\lambda_j} &= 0 \ \text{in}\ M,   \ \ \ u = \text{const.} \ \text{on} \ \  \partial M
\end{align}
and
\begin{align}\label{eq75}
\displaystyle \mathscr{L}_\phi u + au^s\log u + bu^r  &= 0 \ \text{in}\ M,   \ \ \ u  = \text{const.} \ \text{on} \ \  \partial M.
\end{align}

\noindent
{\bf Example 1.} Consider the case
$$H(u)= \sum_{j=1}^r q_ju^{\lambda_j} $$
and compute
$$2H'(u)+(3\varepsilon-2)H(u)/u = \sum_{j=1}^r q_j [2\lambda_j+(3\varepsilon-2)]u^{\lambda_j-1} \le 0$$
provided either $q_j\le 0$ and $\lambda_j\ge 1-3\varepsilon/2$ or $q_j\ge 0$ and $\lambda_j\le 1-3\varepsilon/2$ for each $j$ in the interval $1\le j\le r$.

Note that the quantity $\omega(\varepsilon)=1-3\varepsilon/2$ which serves as bound for $\lambda_j$ can be computed for each $\varepsilon \in (0,1/2]$. For instance, if we allow $\varepsilon \searrow 0$ the quantity $\omega(\varepsilon)$ becomes $\omega(0)=1$ or when picking $\varepsilon = 1/3$ the quantity $\omega(\varepsilon)$ yields $\omega(1/3)=1/2$. By the above discussion we can therefore apply Theorem \ref{thm73} to prove the following result.

\begin{proposition}
Let $M$ be an $n$-dimensional complete smooth metric measure space with compact boundary satisfying $\mathscr{R}ic_\phi^m \ge0$ and $\mathscr{H}_\phi \ge 0$.   Let $u=u(x)$ be any smooth positive bounded solution to \eqref{eq74} such that $u_\nu\ge 0$ over $\partial M$.
Then $u$ must be constant if either $q_j\le 0$ and $\lambda_j\ge 1-3\varepsilon/2$ or $q_j\ge 0$ and $\lambda_j\le 1-3\varepsilon/2$ for each $j$. In particular $\sum_{j=1}^r q_ju^{\lambda_j} \equiv 0$.
\end{proposition}

\noindent
{\bf Example 2.} Consider the case $H(u) = au^s\log u + bu^r$, where $a,b,s$ and $r$ are real. A straightforward calculation yields
$$2H'(u)+(3\varepsilon-2)H(u)/u =a \left[2+[2(s-1)+3\varepsilon]\log u \right]u^{s-1} + b\left[2(r-1)+3\varepsilon\right]u^{r-1}$$
which is nonpositive for $\varepsilon \in (0,1/2]$ if the following assumptions hold:
\begin{align}\label{eq76}
\begin{cases}
a\le 0,  \ \ s \ge 1, \ \ \log u \ge -\frac{2}{2(s-1)+3\varepsilon}\\
\text{or}\\
a\ge 0,  \ \ s \le 1-\frac{3\varepsilon}{2}, \ \ \log u \le -\frac{2}{2(s-1)+3\varepsilon}
\end{cases}
\end{align}
and 
\begin{align}\label{eq77}
b\le 0, \ \ r\ge 1 \ \ \text{or} \ \ b\ge 0, \ \ \ r \le 1- \frac{3\varepsilon}{2}.
\end{align}
Thus an application of Theorem \ref{thm73} yields the following result.
\begin{proposition}
Let $M$ be an $n$-dimensional complete smooth metric measure space with compact boundary satisfying $\mathscr{R}ic_\phi^m \ge0$ and $\mathscr{H}_\phi \ge 0$.   Let $u=u(x)$ be any smooth positive bounded solution to \eqref{eq74} such that $u_\nu\ge 0$ over $\partial M$.
Suppose the assumptions \eqref{eq76} and \eqref{eq77} hold,  then $u$ must be constant.  In particular $au^s\log u + bu^r=0$. Furthermore, if $a\ne 0$ and $b=0$, then $u\equiv 1$.
\end{proposition}

\proof
In consideration of the the assumptions \eqref{eq76} and \eqref{eq77} above we have 
\begin{align*}
\Gamma_+^\varepsilon[u] & = \left[2H'(u)+(3\varepsilon-2)H(u)/u  \right]_+\\
&= \Big[a[2+(2(s-1)+3\varepsilon)\log u ]\Big]_+ \sup_M{u^{s-1}} + \Big[b [2(r-1)+3\varepsilon] \Big]_+\sup_M\{u^{r-1}\}\\
&=0.
\end{align*} 
The claims of the proposition therefore follow at once from Theorem \ref{thm73}.

\qed

 \section*{Declarations}

 \subsection*{Ethics approval and consent to participate}
Not applicable.

\subsection*{Consent for publication}
Not applicable.

\subsection*{Availability of data and material}
 All data generated or analyzed during the study are included in this paper.
 
 \subsection*{Conflict of interests}
 The author declares that there is no competing interests.

\subsection*{Funding}
This work  does not receive any funding.

\subsection*{Author's contributions}
AA is the only contributor in the preparation of this manuscript.

 \subsection*{Acknowledgements} 
This paper was completed during the author research visit to Instituto de Matematicas, Universidad de Granada, Spain.


\begin{thebibliography}{}
%
 		
 \bibitem {[Ab20]}
 		{A. Abolarinwa},
		{\it Gradient estimates for a weighted nonlinear elliptic equation and Liouville type theorems},  J. Geom. Phy. 155 (2020), 103737. 
		
\bibitem{[Ab19]}
			{A. Abolarinwa},  
			{\it Elliptic gradient estimates and Liouville theorems for a weighted nonlinear parabolic equation}, J. Math. Anal. Appl. 473 (2019), 297--312.
			
			
\bibitem{[Ab1]}
		{A. Abolarinwa}, 
			{\it Differential Harnack  and logarithmic Sobolev inequalities along Ricci-harmonic map flow}, Pacific J. Math. 278(2)(2015), 257--290.


\bibitem {[Ab22]}
 		{A. Abolarinwa, A. Ali, F.  Mofarreh},
		{\it Triviality of bounded solutions and gradient estimates for nonlinear $f$-heat equations on complete smooth metric measure spaces},  J. Geom. Phy.  182 (2022) 104670.

\bibitem{[Ab]}{A. Abolarinwa},
		{\it Differential Hanarck estimates for a nonlinear evolution equation of Allen-Cahn 		type}, Mediterr. J. Math. 18, 200(2021).
		
		
 \bibitem {[Ab2]}
			{A. Abolarinwa, S. O. Salawu, C. A. Onate},  
			{\it Gradient estimates for a nonlinear elliptic equation on smooth metric measure spaces and applications}, Heliyon 5(2019).
					
\bibitem{[AE]}
			{A. Abolarinwa, J. O. Ehigie, A. H. Alkhaldi},  
			{\it Harnack inequalities for a class of heat flows with nonlinear reaction 					terms}, J. Geom. Phy. 170 (2021), 104382.
			 			
\bibitem{[Ab0]}
		{A. Abolarinwa,  A. Taheri}, 
		{\it Elliptic gradient estimates for nonlinear $f$-heat equation  on  weighted 				manifolds with time dependent metrics and potentials},  Chaos, Solitons \& Fractals, 142 (2021), 110329.

\bibitem{[AC]}
		S.M. Allen, J.W. Cahn, 
		{\it A microscopic theory for antiphase boundary motion and its application to antiphase domain coarsening},
		 Acta Metall. 27 (1979), 1085-1095.

\bibitem{[AM20]}
		 L. Alzaleq,  V. Manoranjan,  
		{\it  Analysis of Fisher-KPP with a time dependent Allee effect}, 
		 IOP SciNotes 1 (2020) 025003.


 \bibitem {[Ba1]}
			{M. B\v{a}ile\c{s}teanu, X. Cao, A. Pulemotov},  
			{\it Gradient estimates for the heat equation under the Ricci flow}, J. Funct. Anal., 258 (2010), 3517--3542.  
			
		
\bibitem{[BE]}
		{D. Bakry, M. \'Emery}
		{\it Diffusions hypercontractives} In: Az\'ma J., Yor M. (eds) S\'eminaire de Probabilit\'es XIX 1983/84. Lecture Notes in Mathematics, vol 1123. Springer, Berlin, Heidelberg

\bibitem{[BGL]}{D. Bakry, I. Genctil, M. Ledoux}
		{\it Analysis and Geometry of Markov Diffusion Operators}, Springer International Publishing (2014).
		
\bibitem{[Bri]}
K. Brighton,  A Liouville-type theorem for smooth metric measure spaces, J. Geom. Anal. 23 (2013) 562–570.
		
\bibitem{[Cal]}{E.  Calabi}, 
		{\it An extension of Hopf maximum principle with application to Riemannian geometry},  Duke Math. J. 25 (1958), 45--56.
		
\bibitem{[CC]}
		L. Calatroni, P. Colli, 
		{\it Global solution to the Allen-Cahn equation with singular potentials and dynamics boundary conditions}, Nonlinear Anal. 79 (2013), 12--27.
		

\bibitem{[Cao]}{H. D. Cao},
		{\it Recent progress on Ricci solitons in: Recent Advances in Geometric Analysis},
		Adv. Lect. Math. (ALM), 11, International Press, Somerville (2010), 1--38.
		
\bibitem{[CCK]}
	X. Cao,  M. Cerenzia,  D. Kazaras, 
	{\it Harnack estimates for the endangered species equation},
	 Proc. Am. Math. Soc. 143(10) (2014),  4537--4545.

\bibitem{[CLPW]} 
	X. Cao,  B. Liu,  I. Pendleton,  A. Ward, 
	{\it Differential Harnack estimates for Fisher's equation}, 
	 Pac. J. Math. 290(2) (2017),  273--300.

				
\bibitem{[Cas1]}{J. S. Case}, {\it A Yamabe-type problem on smooth metric measure spaces}, J. Differential Geom. 101 (2015), 467--505.

\bibitem{[Cas2]}{J. S. Case, Y-J. Shu, G. Wei}, {\it Rigidity of quasi-Einstein metrics}, Diff. Geom. Appl. 29 (1) (2011), 93--100.

\bibitem{[CaM1]}{D. Castorina, C. Mantegazza},
		{\it Ancient solutions of semilinear heat equations on Riemannian manifolds}, Atti Accad. Naz. Lincei Rend. Lincei Mat. Appl. 28 (1)(2017),  85--101. 

\bibitem{[CaM2]}{D. Castorina, C. Mantegazza},
		{\it Ancient solutions of superlinear heat equations on Riemannian manifolds}, Commun. Contemp. Math. (2020).

\bibitem{[Chen]}{R. Chen},
		{\it Neumann eigenvalue estimate on a compact Riemannian manifold},
		Proc. Amer. Math. Soc. 108 (1990), 961--970.
		
\bibitem{[XChen]}
		X. Chen, 
		{\it Generation and propagation of interface for reaction-diffusion equations}, 
		J. Differ. Equ. 96 (1) (1992),  116 -- 141.

		
\bibitem{[CMJ]} {X.Cheng,T.Mejia, D.T. Zhou},
		{\it Eigenvalue estimate and compactness for closed $f$-minimal surfaces},  Pacific J.  Math. 271 (2014), 347--367.
				
\bibitem{[CZh18]}{Q. Chen, G. Zhao}, 	
		{\it Li-Yau type and Souplet-Zhang type gradient estimates of a parabolic equation for the $V$-Laplacian}, J. Math. Anal. Appl. (2018).	
		
				
\bibitem{[CY75]}{S. Y. Cheng, S. T. Yau}, 
		{\it Differential equations on Riemannian manifolds and their geometric applications}, Comm. Pure Appl. Math. 28 (3)(1975),  333--354. 

\bibitem{[Du23]}
	H.T. Dung, 
	{\it Gradient estimates for a general type of nonlinear parabolic equations under geometric conditions and related problems},  
	Nonlinear Anal. 226 (2023) 113135.
	
			
\bibitem{[Du]}{H. T. Dung},
		{\it Gradient estimates and Harnack inequalities of nonlinear heat equations for the $V$-Laplacian},
		J. Korean Math. Soc. 55(6) (2018), 1285-1305.
								 
\bibitem{[DDW21]}{H. T. Dung, N. T. Dung, J. Y. Wu},
		{\it Sharp gradient estimates on weighted manifolds with compact boundary},
		Commun. Pure Appl. Anal.(2021).
		 
\bibitem{[DK]}{N. T. Dung, N. N. Khanh},
		{\it Gradient estimates for a class of semilinear parabolic equations and their 				applications},
		Vietnam J. Math. (2021).						 
 
 \bibitem{[DKA]}{N. T. Dung, N. N. Khanh, Q. A. Ngo},
		{\it Gradient estimates for some $f$-heat equations driven by Lichnerowicz's equation on complete smooth metric measure spaces},
		Manuscripta Math. 155, (2018), 471--501.	
		
		
\bibitem{[DW21]} {N. T. Dung,  J. Y. Wu}, 
		{\it Gradient estimates for weighted harmonic function with Dirichlet boundary condition},  Nonlinear Anal. 213, (2021), 112498.
		
\bibitem{[Fi]}
		R.A. Fisher, 
		{\it The wave of advance of advantageous genes}, 
		Annu. Eugen. 7 (1937), 355--369.


\bibitem{[FW21]}{X. Fu, J. Y. Wu},
		{\it Gradient estimates for a nonlinear parabolic equation with Dirichlet boundary 	condition},  Kodai  Math.  J.  45 (1)(2022),  96 --109.

\bibitem{[GH]}
	X. Geng,  S. Hou, 
	{\it Gradient estimates for the Fisher-KPP equation on Riemannian manifolds},
	 Bound. Val. Prob.  2018, 25.
	
\bibitem{[Ha]}{ R. Hamilton},  
 {\it The formation of singularities in the Ricci flow}, Surv. Diff. Geom. 2 (1993), 7--136.
 		
 \bibitem{[Ha93]}{ R. Hamilton},  
 {\it A matrix Harnack estimate for the heat equation}, Comm. Anal. Geom. 1(1993), 113--126.
  
  \bibitem{[Hou]} 
 	S. Hou, 
 	{\it Gradient estimates for the Allen-Cahn equation on Riemannian manifolds},
  Proc. Am. Math. Soc. 147,  (2019), 619--628.	
  
 \bibitem{[KPP]}
		A.N. Kolmogorov, I.G. Petrovskii, N.S. Piskunov, 
		{\it A study of the diffusion equation with increase in the amount of substance, and its application to a biological problem},  
		Bull. Moscow Univ. Math. Mech. 1 (1937) 1--26.
  
\bibitem{[KS21]}{K. Kunikawa, Y. Sakurai},
		{\it Yau and Souplet-Zhang type gradient estimates on Riemannian manifolds with boundary under Dirichlet boundary condition},
		Proc.  Amer.  Math.  Soc.  150 (2022),  1767--1777.
		
 \bibitem{[Li]}
	 J. Li,  
	 {\it Gradient estimates and Harnack inequalities for nonlinear parabolic
and nonlinear elliptic equations on Riemannian manifolds},
  J. Funct. Anal. 100 (1991), 233--56.
  
  		
\bibitem{[Li05]}{X. D. Li}, 
{\it Liouville theorems for symmetric diffusion operators on complete Riemannian manifolds}, J. Math. Pures Appl. 84 (2005) 1295--1361
 
\bibitem{[LY86]}{P. Li, S-T. Yau},  
\emph {On the parabolic kernel of the Schr\"odinger operator}, Acta Math. 156 (1986), 153--201.

\bibitem{[LLi1]}{S. Li, X-D. Li},  
\emph {On Harnack inequalities for Witten Laplacian on Riemannian manifolds with super Ricci flows}, Asian J. Math. 22 (2018), 577--598.

\bibitem{[LLi2]}{S. Li, X-D. Li}, 
\emph{Harnack inequalities and $W$-entropy formula for Witten Laplacian on Riemannian manifolds with $K$-super Perelman Ricci flows}, arxiv.org/abs/1412.7034v2

\bibitem{[LLi3]}{ S. Li, X-D. Li}, 
\emph{The $W$-entropy formula for the Witten Laplacian on manifolds with time dependent metrics and potentials}, Pac. J. Math., 278(1), (2015), 173--199

\bibitem{[Lot]}
	J. Lott, 
	{\it Some geometric properties of the Bakry-\'Emery Ricci tensor}, 
	Comment. Math. Helv. 78 (2003) 865--883.
	
\bibitem{[LV]}{J. Lott, C. Villani}, 
	{\it Ricci curvature for metric-measure spaces via optimal transport}, Ann. Math. 169 (2009), 903--991.
	 	
\bibitem{[Ma]}{L. Ma} 
	{\it Gradient estimates for a simple elliptic equation on complete
noncompact Riemannian manifolds}, Journal of Functional Analysis,
241 (1), (2006), 374--382.

\bibitem{[MD10]} {L. Ma,  S.-H. Du}, 
		{\it Extension of Reilly formula with applications to eigenvalue estimates for drifting Laplacians},  
		C. R. Math. Acad. Sci. Paris 348 (2010), 1203--1206.

\bibitem{[MZ]}{B. Ma, F. Zeng}, 
		{\it Hamilton-Souplet-Zhang's gradient estimates and Liouville theorems for a nonlinear parabolic equation}, C. R. Math. Acad. Sci. Paris, Ser. I 356 (5)(2018), 550--557.


\bibitem{[Pe02]}{G. Perelman},
	 {\it The entropy formula for the Ricci Flow and its geometric application}, arXiv:math.DG/0211159v1 (2002).
	 
	 
\bibitem{[Po1]}{P. Pol\'a\v{c}ik, P. Quittner, P. Souplet}, 
		{\it Singularity and decay estimates in superlinear problems via Liouville-type theorems. I. Elliptic equations and systems}, Duke Math. J. 139(3)
(2007), 555--579.


\bibitem{[Po2]}{P. Pol\'a\v{c}ik, P. Quittner, P. Souplet}, 
		{\it Singularity and decay estimates in superlinear problems via Liouville-type theorems. II. Parabolic equations}, Indiana Univ. Math. J. 56(2) (2007),  879-908.
		

\bibitem{[Olive]}{X. Ramos Oliv\'e},
		{\it Neumann Li-Yau gradient estimate under integral Ricci curvature bounds},
		Proc. Amer. Math. Soc. 147 (2019), 411--426.


\bibitem{[Rei]} {R.C. Reilly}, 
		{\it Applications of the Hessian operator in a Riemannian manifold}, 
		Indiana Univ. Math. J. 26 (1977), 459--472.		
		
\bibitem{[Sak2]}{Y. Sakurai},
		{\it Rigidity of manifolds with boundary under a lower Bakry-\'Emery Ricci curvature 		bound},
		Tohoku Math. J. 71 (2019), 69--109.
		
\bibitem{[SY94]}{R. Schoen,  S.T. Yau}, 
		{\it  Lectures on Differential Geometry},  International Press, 1994.

\bibitem{[SZ06]}{P. Souplet, Q. S. Zhang},
		{\it Sharp gradient estimate and Yau's Liouville theorem for the heat equation on noncompact manifolds}, Bull. London Math. Soc. 38 (6)(2006), 1045--1053.
		
				
\bibitem{[Vil]}{C. Villani},
	 {\it Optimal Transport: Old and New, A Series of Comprehensive studies in mathematics}, vol. 338. Springer, Berlin

\bibitem{[Wang97]}{J. P. Wang},
		{\it Global heat kernel estimates},
		Pacific J. Math. 178 (1997), 377--398.

\bibitem{[WZZ16]}{L.-F. Wang, Z.-Y. Zhang, Y.-J.Zhou},
		{\it Comparison theorems on smooth metric measure spaces with boundary},
		Adv. Geom. 16 (2016), 401--411.

		 	  
\bibitem{[WeW09]}{G. F. Wei, W. Wylie},
\emph{Comparison geometry for the Bakry-\'Emery Ricci tensor}, J. Differ. Geom. 83 (2009), 377--405.

\bibitem{[Wu15]}{J-Y. Wu},
	 {\it Elliptic gradient estimates for a weighted heat equation and applications}, Math. Z.(2015).
	 
\bibitem{[Wu]}{J-Y. Wu},
	 {\it Gradient estimates for a nonlinear parabolic equation and Liouville theorems},
	 Manuscript Math. 159 (2018), 511--547.
	 
\bibitem{[Ya]}
	Y. Yang,
	{\it Gradient estimates for a nonlinear parabolic equation on Riemannian manifolds},
	 Proc. Am. Math. Soc. 136 (2008),  4095 -- 4102.

  \bibitem{[Ya1]} {S-T. Yau}, 
		{\it Harmonic functions on complete  Riemannian manifolds}, 
Comm. Pure Appl. Math. 28 (1975), 201--228.
 
\bibitem{[Ya2]} {S-T. Yau}, 
		{\it Some function-theoretic properties of complete Riemannian manifolds and their applications to geometry}, Indiana Math. J. 25 (1976), 659--670.
 
 \end{thebibliography}
\end{document}